\documentclass[11pt]{article}
\usepackage{amsmath, amssymb, amsthm, amsbsy, ifthen, graphicx}
\date{}
\title{Ramsey games with giants}
\author{
Tom Bohman \thanks{Department of Mathematical Sciences, Carnegie
Mellon University, Pittsburgh, PA 15213, email:
tbohman@math.cmu.edu. Research supported in part by NSF award
DMS-0701183.} 
\and 
Alan Frieze \thanks{Department of Mathematical
Sciences, Carnegie Mellon University, Pittsburgh, PA 15213, email:
alan@random.math.cmu.edu. Research supported in part by NSF award
DMS-0753472.} 
\and 
Michael Krivelevich \thanks{School of
Mathematical Sciences, Raymond and Beverly Sackler Faculty of Exact
Sciences, Tel Aviv University, Tel Aviv 69978, Israel, e-mail:
krivelev@post.tau.ac.il. Research supported in part by USA-Israel
BSF grant 2006322, by grant 1063/08 from the Israel Science
Foundation, and by a Pazy memorial award.} 
\and 
Po-Shen Loh
\thanks{Department of Mathematics, Princeton University, Princeton,
NJ 08544, e-mail: ploh@math.princeton.edu. Research supported in
part by a Fannie and John Hertz Foundation Fellowship, an NSF
Graduate Research Fellowship, and a Princeton Centennial
Fellowship.} 
\and 
Benny Sudakov \thanks{Department of Mathematics,
UCLA, Los Angeles, CA 90095, email: bsudakov@math.ucla.edu. Research supported in part by
NSF CAREER award DMS-0812005 and by a USA-Israeli BSF grant.} }

\newtheorem{theorem}{Theorem}[section]
\newtheorem{definition}[theorem]{Definition}
\newtheorem{fact}[theorem]{Fact}
\newtheorem{proposition}[theorem]{Proposition}
\newtheorem{lemma}[theorem]{Lemma}
\newtheorem{corollary}[theorem]{Corollary}

\newtheorem*{customtheorem}{Theorem}

\newcommand{\pr}[1]{\mathbb{P}\left[#1\right]}
\newcommand{\E}[1]{\mathbb{E}\left[#1\right]}
\newcommand{\bin}{\text{Bin}}

\newcommand{\whp}{\textbf{whp}}

\newcommand{\orient}{\textsc{orient}}
\newcommand{\Gnpd}{\overrightarrow{G}_{n,p}}
\newcommand{\Gd}{\overrightarrow{G}}
\newcommand{\Ft}{e_1, \ldots, e_t}
\newcommand{\alignsp}{\quad\quad\quad\quad}

\def\cE{{\cal E}}
\begin{document}
\maketitle

\begin{abstract}
  The classical result in the theory of random graphs,
  proved by Erd\H{o}s and R\'enyi in 1960, concerns the threshold for
  the appearance of the giant component in the random graph process.
  We consider a variant of this problem, with a Ramsey flavor.
  Now, each random edge
  that arrives in the sequence of rounds must be colored with one of
  $r$ colors.  The goal can be either to create a giant component in
  every color class, or alternatively, to avoid it in every color.  One
  can analyze the offline or online setting for this problem.
  In this paper, we consider all these variants and provide 
  nontrivial upper and lower bounds; in certain cases (like online
  avoidance) the obtained bounds are asymptotically tight.
\end{abstract}

\section{Introduction}

Let $G_{n,m}$ be the Erd\H{o}s-R\'enyi random graph with $n$ labeled
vertices and $m$ randomly chosen edges. A celebrated result of
Erd\H{o}s and R\'enyi, probably the single most important result in
the theory of random graphs, discovered a threshold for the
appearance of the giant component in this random model. Erd\H{o}s
and R\'enyi proved that if $m\le (1-\epsilon)\frac{n}{2}$ for a constant
$\epsilon>0$, then \whp\footnote{As customary, we write that a graph
property $\mathcal{P}$ holds with high probability, or \whp\ for brevity, 
if the probability of $\mathcal{P}$ tends to
1 as the number of vertices $n$ tends to infinity.} the random graph
$G_{n,m}$ has all of its connected components of order at most
logarithmic in $n$; on the other hand, if $m\ge (1+\epsilon)\frac{n}{2}$
then \whp\ $G_{n,m}$ has a unique connected component of linear
size, the so called giant component, while all other components are
at most logarithmic in size. This result can be formulated
equivalently in terms of the random graph process: if the process
starts with the empty graph $G_0$ on $n$ vertices, and at stage
$i\ge 1$ a random missing edge is added to $G_{i-1}$ to form $G_i$,
then after the first $(1-\epsilon)\frac{n}{2}$ rounds the resulting graph
typically has all connected components of at most logarithmic size,
while after $(1+\epsilon)\frac{n}{2}$ rounds \whp\ the unique giant
component is born, while all other components are of size $O(\log
n)$. Since then, there have been numerous extensions to this 
fundamental result. One further ramification
is considered in this paper.

Recently, quite a lot of attention and research effort has been
devoted to controlled random graph processes. In processes of this
type, an input graph or a graph process is usually generated fully
randomly, but then an algorithm has access to this random input
and can manipulate it in some well defined way (say, by dropping
some of the input edges, or by coloring them), aiming to achieve
some preset goal. There is usually the so called {\em online}
version where the algorithm must decide on its course of action
based only on the history of the process so far and without assuming
any familiarity with future random edges, and the {\em offline}
version, where the algorithm has access to the whole history of
the process and makes its decisions based on the full knowledge
of the process. We will give corresponding accurate definitions for
our setting later.

Applied to the question about the appearance of the giant component,
the first such version chronologically is probably the so-called
\emph{Achlioptas process}. This process is named after Dimitris Achlioptas,
who posed the following question about 10 years ago. Suppose random edges
arrive in pairs, and an online algorithm can choose one of them, put
it into the graph, and return the other edge to the pool. Is it
possible to design an algorithm that \whp\ delays the appearance of
the giant components for noticeably longer than the Erd\H{o}s-R\'enyi
$0.5n$ steps? This question was answered affirmatively in \cite{BF} by the
first two authors of the present paper, who exhibited an
algorithm that \whp\ survives for at least $0.535n$ rounds without
creating the giant component. Since then, there has been a series of
papers about the giant component in Achlioptas
processes, where a variety of scenarios and goals (online and
offline algorithms, delaying or accelerating the appearance of the
giant component) have been considered.

Here we consider a Ramsey-type version of controlled random
processes. In this version, incoming random edges are colored by an
algorithm in one of $r$ colors, for a fixed $r\ge 2$.  The goal of
the algorithm is to achieve or maintain a certain monotone graph property in
all of the colors. 
This setting originates in the papers
of R\"odl and Ruci\'nski \cite{RoRu0, RoRu}, who determined when
$G_{n,m}$ satisfies the Ramsey property of having a monochromatic
copy of a fixed graph $H$ in any $r$-coloring of the edges. In our
terminology they considered the offline version of the problem, and
the property $\mathcal{P}$ to avoid in each color was the appearance of a copy
of a fixed graph $H$. The online version of the problem for the case
of two colors and $H=K_3$ was treated by Friedgut, Kohayakawa,
R\"odl, Ruci\'nski and Tetali in \cite{FKRRT}, and extended to
a wider variety of graphs by Marciniszyn,
Sp\"ohel and Steger in \cite{MSS1, MSS2}. The online setting of
achieving Hamiltonicity in each of $r$ colors has been addressed in
\cite{KLuS}.

In the present paper, we investigate several Ramsey-type problems
involving the giant component.  We consider whether or not it is
possible to color the edges of $G_{n,m}$ in $r$ colors with the
objective of creating a giant component in every color class, or of
avoiding a giant component in every color.  We study both the
offline and online settings. In the offline setting, an algorithm
gets access to the entire graph, generated according to the
probability distribution $G_{n,m}$; in the online setting the edges
of $G_{n,m}$ are first ordered in a random order and then revealed
to the algorithm one by one (i.e., the algorithm observes the random
graph process and colors each new edge as it arrives).

The main objective of this paper is to show new interesting
questions, and not necessarily to get precise answers to all of
them.  We do determine the offline thresholds for these problems for
all values of $r$, but the online setting remains open.  There, we
show that for two colors, there is always a separation phenomenon
away from the trivial bounds, and then calculate asymptotic bounds
for large numbers of colors.

As a warm-up, consider the offline threshold for creating a giant in
every color.  Recall that if $m < (1-\epsilon) \frac{n}{2}$ for any
fixed $\epsilon > 0$, then \whp\ $G_{n,m}$ itself has all components
of size $O(\log n)$.  On the other hand, one can show that for
$m>(1+\epsilon) \frac{n}{2}$, \whp\ it is possible to color the
edges of $G_{n,m}$ with any fixed number of colors $r \geq 2$, so
that every color class contains a component of order $\Omega(n)$.
Indeed, Ajtai, Koml\'os, and Szemer\'edi 
proved in \cite{AKS} that \whp\ $G_{n,(1+\epsilon)\frac{n}{2}}$ contains a path of
length $c_\epsilon n$.  (Here and later in the paper, we will write
$c_\epsilon$ to specify a positive constant determined only by
$\epsilon$.)  By splitting this path into $r$ paths of length
$c_\epsilon n/r$, the result follows.

The question of avoiding giants in all colors offline is not so
simple. It turns out that the threshold for avoiding giants in $r$
colors is precisely the same as that of \emph{$r$-orientability},
which says that it is possible to direct all of the edges of the
graph so that the resulting digraph has maximum in-degree at most
$r$.  Cain, Sanders and Wormald \cite{CSW}, and Fernholz and
Ramachandran \cite{FR} recently discovered that this threshold
coincides with the number of edges needed to make the $(r+1)$-core
have average degree above $2r$.  More precisely, they showed that
for any integer $r \geq 2$, there is an explicit threshold $\psi_r$
such that the following holds.  For any $\epsilon > 0$, if $m >
(\psi_r + \epsilon)n$, then \whp\ $G_{n,m}$ contains a subgraph with
average degree at least $2r + c_\epsilon$, where $c_\epsilon > 0$.
On the other hand, if $m < (\psi_r - \epsilon)n$, then $G_{n,m}$ is
$r$-orientable \whp.  The asymptotic dependence of $\psi_r$ on $r$
is $\psi_r = r - \frac{1}{2}\big(\frac{2}{e} + o(1) \big)^r$, as
calculated in \cite{CSW}.  We now state our first main theorem in
terms of this threshold.

\begin{theorem}
  \label{thm:offline-avoid}
  Given any fixed $r$, let $\psi_r$ be the threshold referenced above.
  For any $\epsilon > 0$, if $m < (\psi_r - \epsilon)n$, then \whp\ it
  is possible to color the edges of $G_{n,m}$ with $r$ colors such
  that each color class contains components of order only $o(n)$.  On
  the other hand, if $m > (\psi_r + \epsilon)n$, then \whp\ every
  $r$-edge-coloring of $G_{n,m}$ has a color class with a component of
  order at least $c_\epsilon n$.
\end{theorem}

\noindent \textbf{Remark.}\, This was also recently and independently
discovered by Sp\"ohel, Steger, and Thomas \cite{SST}.

\vspace{3mm}

We also consider online versions of these problems, in which the $m$
edges come sequentially, and each must be colored as soon as it
appears.  Precisely, we consider the process to be a sequence of $m$
rounds.  In each round, a random edge arrives, independently and
uniformly distributed over all pairs of vertices.  If it repeats an
existing edge, then we do not force ourselves to recolor it.  This
is not an important issue, because we will never consider more than
$O(n)$ rounds, but it is more convenient to use this product
probability space with full independence between the rounds.

Here, we have several results.  First we state them for avoiding
giants in all colors.  The offline upper bound of course supplies an
upper bound for the online case as well.  Indeed, a standard
coupling argument (Fact \ref{fact:coupling} in the next section)
translates the offline upper bound to the case where the rounds have
independent edges (possibly with repetitions).  So, after $(\psi_r +
\epsilon) n$ rounds, \whp\ every possible coloring of them contains
a giant component, where the dependence of $\psi_r$ on $r$ is
$\psi_r = r - \frac{1}{2}\big(\frac{2}{e} + o(1) \big)^r$.

On the other hand, by taking the natural online adaptation of the
offline avoidance strategy, which was based on edge orientation, we
found a randomized online algorithm which matches the first-order
asymptotic of $\psi_r = (1-o(1))r$.

\begin{theorem}
  \label{thm:online-avoid}
  For any $\epsilon > 0$, the following holds for all sufficiently
  large $r$.  There is an online randomized algorithm which can last
  for $(1 - \epsilon) rn$ rounds, while keeping all connected
  components in each of $r$ color classes smaller than $o(n)$ \whp.
\end{theorem}

For large $r$, this is asymptotically a factor of 2 better than the
trivial bound of $(1-\epsilon) \frac{rn}{2}$ rounds, obtained by
coloring each edge independently at random.  On the other hand, for
small $r$, beating the trivial bound corresponds to using $\epsilon =
0.4999$ above (say), and this requires $r > 50$.  For the extreme case
of small $r$, we have the following result using an entirely different
strategy, which improves upon the trivial bound for all $r$ by a
factor of approximately $1.06$.

\begin{theorem}
  \label{thm:online-avoid-2}
  There is an online algorithm which can 2-color edges for $1.06 n$
  rounds, while keeping all connected components in both color classes
  of size at most $O(\log n)$ \whp.
\end{theorem}

\noindent \textbf{Remark.}\, Although the theorem is stated only for
$r=2$, it immediately gives a strategy for all even $r$, by splitting
the colors into $\frac{r}{2}$ pairs.  At each round, one of the color
pairs is randomly chosen, and the above algorithm is used to decide
which of the two colors in the pair to use.  Then, this will avoid
giants in all colors for $1.06n \cdot \frac{r}{2}$ rounds \whp.  For
odd $r$, one can run the above modification for $1.06n \cdot
\frac{r-1}{2}$ rounds using only the first $r-1$ colors, and then an
additional $(1-\epsilon) \frac{n}{2}$ rounds using only the $r$-th
color.  This beats the trivial bound of $\frac{rn}{2}$ by a factor
which approaches $1.06$ as $r$ grows.

\vspace{3mm}

When the objective is to create giants in every color class, the
trivial bounds are as follows.  Certainly, if fewer than
$(1-\epsilon)\frac{n}{2}$ edges are observed, then \whp\ there will be no
giant in the uncolored graph, so one cannot hope to create $r$
monochromatic giants any faster. Note that this trivial lower bound
turned out to be the truth in the offline setting, even though it
does not grow with $r$.  We will show that in the online case, there
is a lower bound which does.

\begin{theorem}
  \label{thm:online-embrace-lower}
  There is a constant $c \approx 0.043$ such that after $(c \log_2 r)
  n$ edges are $r$-colored by any online algorithm, \whp\ some color
  class still has all components of order only $O(\log n)$.  For
  $r=2$, the same result holds for $c' n$ edges for any $c' <
  2-\sqrt{2} \approx 0.586$.
\end{theorem}

On the other hand, the trivial strategy of randomly coloring each
edge succeeds when the number of edges surpasses $rn/2$.  We are
able to give an online algorithm which asymptotically performs far
better than the trivial one.

\begin{theorem}
  \label{thm:online-embrace-upper}
  There is an online algorithm such that for any $\epsilon > 0$,
  after $(c_r + \epsilon) n$ edges every color class contains a
  connected component of order at least $c_\epsilon n$ \whp, where the
  dependence of $c_r$ on $r$ is $c_r = (1+o(1))\frac{\sqrt{r}}{2}$.
\end{theorem}

For the specific case of 2 colors, one can adapt the argument and
obtain a value of $c_2 = \frac{3}{4}$, but we give a slightly more
sophisticated strategy which creates giants even faster.

\begin{theorem}
  \label{thm:online-embrace-upper-2}
  There is an online algorithm such that for any $\epsilon > 0$,
  after $0.733 n$ rounds both color classes
  contain connected components of order at least $c_\epsilon n$
  \whp.
\end{theorem}

This paper is organized as follows.  The next section reviews some
standard probabilistic facts, and then develops a general tool which
extends a recent result of Spencer and Wormald from \cite{SW}.  This
allows us to control the evolution of the susceptibility of a graph
under the addition of random edges.  Section \ref{sec:offline-avoid}
completely resolves the offline case, by proving Theorem
\ref{thm:offline-avoid}.  For the online setting, Sections
\ref{sec:online-avoid} and \ref{sec:online-create} consider the
respective problems of avoiding and creating giants.  The final
section contains some concluding remarks.

Throughout our paper, we will omit floor and ceiling signs whenever
they are not essential, to improve clarity of presentation.  All
logarithms are in base $e \approx 2.718$ unless otherwise specified.
The following asymptotic notation will be utilized extensively.  For
two functions $f(n)$ and $g(n)$, we write $f(n) \ll g(n)$, $f(n) =
o(g(n))$, or $g(n) = \omega(f(n))$ if $\lim_{n \rightarrow \infty}
f(n)/g(n) = 0$, and $f(n) = O(g(n))$ or $g(n) = \Omega(f(n))$ if
there exists a constant $M$ such that $|f(n)| \leq M|g(n)|$ for all
sufficiently large $n$. The number of vertices $n$ is assumed to be
sufficiently large where necessary.

%

\section{Preliminaries}

In this section, we review some standard facts commonly used in
Probabilistic Combinatorics.  Then, we use them to prove a useful
result (Theorem \ref{thm:evolve-suscep}) which shows that a certain
graph parameter, the \emph{susceptibility}, tracks a natural
differential equation.  This extends a result of Spencer and
Wormald, and we state it in a general-purpose form for the
convenience of possible future citations.

\subsection{Probabilistic tools}

We recall the Chernoff bound for exponential concentration of the
binomial distribution.  The following formulation appears in, e.g.,
\cite{AS}.

\begin{fact}
  \label{fact:conc}
  For any $\epsilon > 0$, there exists $c_\epsilon > 0$ such that any
  binomial random variable $X$ with mean $\mu$ satisfies $\pr{| X -
    \mu | > \epsilon \mu} < e^{-c_\epsilon \mu}$.
\end{fact}

A binomial random variable is the sum of independent indicator
variables.  We also need concentration in settings without complete
independence.  Recall that a martingale is a sequence $X_0, X_1,
\ldots$ of random variables such that each conditional expectation
$\E{X_{t+1} \mid X_0, \ldots, X_t}$ is precisely $X_t$.  The
Hoeffding-Azuma inequality (see, e.g., \cite{AS}) provides
concentration for martingales with bounded step-wise increments
$|X_{t+1} - X_t|$, and this has been widely used in probabilistic
combinatorics.

When only one-sided concentration is needed, it can be convenient to
consider instead a \emph{supermartingale}, which only requires
$\E{X_{t+1} \mid X_0, \ldots, X_t} \leq X_t$ for all $t$.  We will use
the analogue of Hoeffding-Azuma for supermartingales, which follows
from exactly the same proof as for martingales (see, e.g., \cite{M} or
\cite{W95}).

\begin{fact}
  \label{fact:azuma}
  Let $X_0, \ldots, X_n$ be a supermartingale, with bounded
  differences $|X_{i+1} - X_i| \leq C$.  Then for any $\lambda \geq
  0$,
  \begin{displaymath}
    \pr{X_n \geq X_0 + \lambda}
    \ \leq \
    \exp\left\{-\frac{\lambda^2}{2 C^2 n}\right\}.
  \end{displaymath}
\end{fact}

We can also define \emph{submartingales}\/ via the requirements
$\E{X_{t+1} \mid X_0, \ldots, X_t} \geq X_t$; estimates on their lower
tails, similar to the above fact, follow by symmetry.

Finally, we will frequently switch between the models $G_{n,p}$,
$G_{n,m}$, and the product space of $m$ independent uniform random
edges, depending on which one is the most convenient.  Adding more
edges makes it harder to avoid giants, but easier to create them, so
all properties we consider are monotone.  Hence the following fact
allows us to translate results between the models, while still
keeping everything sharp to first-order.

\begin{fact}
  \label{fact:coupling}
  Fix any constant $\epsilon > 0$, and suppose that $m = m(n)$ tends
  to infinity with $n$, but $m = o(n^2)$.  Then there are couplings of
  the corresponding probability spaces such that the following hold.
  \begin{description}
  \item[(i)] $G_{n,m} \subset G_{n,p}$ \whp\ for $p =
    (1+\epsilon)\frac{2m}{n}$, and $G_{n,m} \supset G_{n,p}$ \whp\ for
    $p = (1-\epsilon)\frac{2m}{n}$.
  \item[(ii)] The graph formed by generating $m$ random edges
    (possibly with repetition) is always contained in $G_{n,m}$, and
    \whp\ contains $G_{n,m'}$ with $m' = (1-\epsilon)m$.
  \end{description}
\end{fact}

\noindent \textbf{Proof sketch.}\, By the standard coupling of
$G_{n,m}$ and $G_{n,p}$ via the random graph process, part (i) follows
from the Chernoff bound on $\bin\big[\binom{n}{2}, p\big]$.  For part
(ii), one can similarly couple $G_{n,m'}$ with the product space of
$m$ edges by considering an infinite sequence of independent random
edges.  Then, the $m$-edge product space is the projection onto the
first $m$ choices, and $G_{n,m'}$ is the graph consisting of the first
$m'$ distinct edges.  So, it suffices to show that \whp, there are at
least $(1-\epsilon) m$ distinct edges among the first $m$ sampled with
replacement.  Observe that when the $k$-th edge is sampled, the
probability that it is a repetition of a previously sampled edge is
always less than $k/ \binom{n}{2} < \frac{\epsilon}{2}$ since $m =
o(n^2)$.  Therefore, the number of samples which are repetitions is
stochastically dominated by $\bin\big[m, \frac{\epsilon}{2} \big]$,
which is at most $\epsilon m$ \whp\ by the Chernoff bound.  Then, the
number of distinct edges is at least $(1-\epsilon)m$, as desired.
\hfill $\Box$

\subsection{Evolution of susceptibility}

One of the most useful parameters for studying the giant component
of a graph is the \emph{susceptibility}\/.  For a graph $G$, this is
defined as $S(G) = \frac{1}{n} \sum_v C_v$, where $C_v$ is the size
of the connected component in $G$ containing $v$.  Note that this
also equals $\frac{1}{n}$ times the sum of the squares of the
component sizes.  Many researchers have investigated the evolution
of the susceptibility under random edge addition, starting with
Bohman and Kravitz, who used this to analyze the Achlioptas process
in \cite{BKr}.

More recently, Spencer and Wormald proved in \cite{SW} that for $m$ up
to $(1-\epsilon)\frac{n}{2}$, the susceptibility of the $m$-edge
random graph evolves like the solution $\phi(m)$ of the differential
equation $\phi' = \frac{2}{n} \phi^2$ with initial condition $\phi(0)
= 1$.  The heuristic for this differential equation is quite natural,
although the formal proof is nontrivial.  Indeed, when a random edge
is added to some intermediate (and subcritical) $G$, its endpoints
typically lie in different components, each of which has expected size
$S(G)$.  If both component sizes are close to $S(G)$, then the
increment to $S(G)$ after adding the edge is roughly $\frac{1}{n}
\big[ (S(G)+S(G))^2 - 2S(G)^2 \big] = \frac{2}{n} S(G)^2$.  Thus, one
might expect the evolution of $S(G)$ to follow $\phi' = \frac{2}{n}
\phi^2$.  The solution of this differential equation is $\phi(m) =
\big(1 - \frac{2}{n} m\big)^{-1}$, so it only blows up when $m$
reaches $\frac{n}{2}$.  This matches the classical threshold of the
giant component, because the result of Spencer and Wormald
concentrates $S(G_{n,m})$ around $\phi(m)$ for $m$ up to
$(1-\epsilon)\frac{n}{2}$.  In this range, $S(G_{n,m})$ is then
bounded by a constant, and we can always trivially bound the size of
the largest component by $\sqrt{n S(G)}$, so the largest component is
$o(n)$ \whp.

However, once we start to color edges, the color classes are no longer
Erd\H{o}s-R\'enyi random graphs.  It is then crucial to control the
evolution of susceptibility from initial graphs which are non-empty.
One of the main contributions of \cite{SW} was a result of this
nature, but it only controlled one phase of evolution.  In order to
formulate it, we need the following definition.

\begin{definition}
  A graph has a $\boldsymbol{K, c}$ \textbf{component tail} if for all
  positive integers $s$, at most $Ke^{-cs}$-fraction of its vertices lie
  in components of order at least $s$.
\end{definition}

Note that a $K,c$ component tail immediately implies that all
components have order $O(\log n)$.  Now we restate a key result of
Spencer and Wormald (Theorem 3.1 of \cite{SW}), translated into an
equivalent form via Fact \ref{fact:coupling}.

\begin{fact}
  \label{fact:SW}
  Let $L, K, c, \gamma$ be positive real numbers.  Let $G$ be a graph
  on $n$ vertices with a $K,c$ component tail and $S(G) \leq L$.  Add
  $(1-\gamma) \frac{n}{2L}$ independent random edges to $G$, ignoring
  repeated edges, and let the result be $G'$.  Then there exist $K',
  c'$ such that $G'$ has a $K', c'$ component tail \whp.
\end{fact}

The $K', c'$ component tail is very useful, because it bounds the
entire distribution of the component sizes.  However, our arguments
also need control of the new value of the susceptibility after random
edge addition, so we prove the following extension of the above
result. This can be done using the methods used in \cite{SW}, but we
include here an alternate (and simpler) proof, following ideas from
\cite{B}.

\begin{theorem}
  \label{thm:evolve-suscep}
  Let $L, K, c, \gamma$ be positive real numbers.  Let $G$ be a graph
  on $n$ vertices with a $K,c$ component tail and $S(G) \leq L$.  Add
  $(1-\gamma) \frac{n}{2L}$ independent random edges to $G$, ignoring
  repeated edges, and let the result be $G'$.  Then there exist $K',
  c'$ such that \whp\ $G'$ has a $K', c'$ component tail, \textbf{and
    $\boldsymbol{S(G') \leq \frac{L}{\gamma} + o(1)}$}.
\end{theorem}

\noindent \textbf{Remark.}\, The bound $\frac{L}{\gamma}$ arises
from the following heuristic.  Suppose that the initial
susceptibility is $L$.  We will show that its evolution is dictated
by the differential equation $\phi' = \frac{2}{n} \phi^2$ with
initial condition $\phi(0) = L$, whose solution is $\phi(t) =
\big(\frac{1}{L} - \frac{2}{n} t\big)^{-1}$. Substituting $t =
(1-\gamma) \frac{n}{2L}$ gives $\frac{L}{\gamma}$.

\vspace{3mm}

\noindent \textbf{Proof.}\, Note that by definition, the
susceptibility is always at least 1, so we will implicitly use $L \geq
1$ throughout the proof.  Let $T = (1-\gamma) \frac{n}{2L}$.  Let
$e_1, \ldots, e_T$ denote the incoming edges, and let $G_t$ be the
graph after the addition of the first $t$ of them.  Fact \ref{fact:SW}
gives constants $K', c'$ such that $G_T$ has a $K', c'$ component tail
\whp.

Let $\phi(t) = \big(\frac{1}{L} - \frac{2}{n} t\big)^{-1}$.  We now
formalize our heuristic argument which suggests that $S(G_t)$ evolves
like $\phi(t)$.  For each $t$, let $\mathcal{E}_t$ be the event that $G_t$ has a
$K', c'$ component tail and $S(G_t) \leq \phi(t) +
e^{\frac{5L}{\gamma} \frac{t}{n}} n^{-\frac{1}{3}}$.  Note that we
will only run $t$ up to $T \leq n$, so the exponential factor is only
at most a constant, and hence the error term tends to zero as $n$
grows.  Now, consider the sequence of random variables:
\begin{displaymath}
  X_t \ = \ \begin{cases}
    S(G_t) - \phi(t) - e^{\frac{5L}{\gamma} \frac{t}{n}} n^{-\frac{1}{3}} & \text{ if $\mathcal{E}_{t-1}$ holds,} \\
    X_{t-1} & \text{ otherwise. }
  \end{cases}
\end{displaymath}
We claim that $X_t$ is a supermartingale.  Indeed, suppose that $G_t$
has components of order $C_1, C_2, \ldots$ If the incoming edge $v_1
v_2$ has $v_1$ in the $i$-th component and $v_2$ in the $j$-th
component, then the susceptibility increases by exactly $\frac{1}{n} [
(C_i + C_j)^2 - C_i^2 - C_j^2 ] = \frac{2}{n} C_i C_j$ when $i \neq
j$, and zero otherwise.  Therefore,
\begin{eqnarray*}
  \E{S(G_{t+1}) \mid \Ft}
  &=& S(G_t) + \sum_{i \neq j} \frac{2}{n} C_i C_j \cdot \frac{C_i}{n} \frac{C_j}{n-1} \\
  &\leq& S(G_t) + \frac{2}{n-1} \left( \frac{1}{n} \sum_i C_i^2 \right)^2 \\
  &=& S(G_t) + \frac{2}{n-1} S(G_t)^2.
\end{eqnarray*}
We use this to bound the expected conditional increment in $X_t$.
Note that for the purposes of bounding $\E{X_{t+1} \mid \Ft}$ we may
assume that $\mathcal{E}_t$ holds (otherwise this conditional expectation is
trivially equal to $X_t$).  Using the above, and the convexity of
$\phi$ and the exponential, we have:
\begin{align*}
  &\E{X_{t+1} - X_t \mid \Ft, \mathcal{E}_t} \\
  &\alignsp \leq \ \frac{2}{n-1} S(G_t)^2 - \left(\phi(t+1) - \phi(t)\right) -
  \left( e^{\frac{5L}{\gamma} \frac{t+1}{n}} - e^{\frac{5L}{\gamma} \frac{t}{n}} \right) n^{-\frac{1}{3}} \\
  &\alignsp \leq \ \frac{2}{n-1} S(G_t)^2 - \phi'(t) -
  \frac{5L}{\gamma} \frac{1}{n} e^{\frac{5L}{\gamma} \frac{t}{n}} n^{-\frac{1}{3}} \\
  &\alignsp = \ \frac{2}{n-1} S(G_t)^2 - \frac{2}{n} \phi(t)^2 -
  \frac{5L}{\gamma} \frac{1}{n^{4/3}} e^{\frac{5L}{\gamma} \frac{t}{n}} \\
  &\alignsp \leq \ \frac{2}{n-1} \left( \phi(t) + e^{\frac{5L}{\gamma} \frac{t}{n}} n^{-\frac{1}{3}} \right)^2
  - \frac{2}{n} \phi(t)^2 -
  \frac{5L}{\gamma} \frac{1}{n^{4/3}} e^{\frac{5L}{\gamma} \frac{t}{n}} \\
  &\alignsp = \ \frac{2}{n(n-1)} \phi(t)^2 + \frac{4}{(n-1)n^{1/3}} \phi(t) e^{\frac{5L}{\gamma} \frac{t}{n}}  + \frac{2}{(n-1)n^{2/3}} e^{\frac{10L}{\gamma} \frac{t}{n}} - \frac{5L}{\gamma} \frac{1}{n^{4/3}} e^{\frac{5L}{\gamma} \frac{t}{n}}.
\end{align*}
We will only run $t$ up to $T = (1-\gamma) \frac{n}{2L}$, so we always
have $\frac{t}{n} < 1$, as well as $\phi(t) \leq \frac{L}{\gamma}$
because $\phi$ is increasing.  Plugging in these bounds, the $\phi(t)$
and exponential factors are replaced by constants, so the asymptotic
behavior of each term is determined by the power of $n$ in the
denominator.  Hence the second and fourth terms dominate, giving
\begin{eqnarray*}
  \E{X_{t+1} - X_t \mid \Ft, \mathcal{E}_t}
  &\leq& (1+o(1)) \left( \frac{4}{n^{4/3}} \frac{L}{\gamma} e^{\frac{5L}{\gamma}}
    - \frac{5L}{\gamma} \frac{1}{n^{4/3}} e^{\frac{5L}{\gamma}} \right) \\
  &=& - (1+o(1)) \frac{L}{\gamma n^{4/3}} e^{\frac{5L}{\gamma}},
\end{eqnarray*}
which is negative for sufficiently large $n$.  Therefore, $X_t$ is
indeed a supermartingale.  Observe that $X_0 = -n^{-1/3}$.  We will
use the Hoeffding-Azuma inequality (Fact \ref{fact:azuma}) to prove
that \whp, $X_t < 0$ for every $t \leq T$.  For this, note that the
one-step change in $X_t$ is zero if $G_t$ does not have a $K', c'$
component tail.  Otherwise, as previously remarked, all components of
$G_t$ are bounded by some $C \log n$, so the maximum change in the
susceptibility is $\frac{2}{n} (C \log n)^2$.  To bound the one-step
change in the error term $\phi(t) + e^{\frac{5L}{\gamma} \frac{t}{n}}
n^{-\frac{1}{3}}$, which is an increasing convex function, it suffices
to use the first derivative at $t=T$.  Recalling that $T =
(1-\gamma)\frac{n}{2L}$, this turns out to be precisely
\begin{displaymath}
  \left. \frac{d}{dt} \right|_{t=T}
  \ = \
  \left[ \left( \frac{1}{L} - \frac{2}{n} T \right)^{-2} \cdot \frac{2}{n} \right] +
  \left[ e^{\frac{5L}{\gamma} \frac{T}{n}} n^{-1/3} \cdot \frac{5L}{\gamma n} \right]
  \ = \
  \left( \frac{\gamma}{L} \right)^{-2} \frac{2}{n}   +
  e^{\frac{5L}{\gamma} \frac{T}{n}} \frac{5L}{\gamma n^{4/3}},
\end{displaymath}
which is clearly $O(n^{-1})$ because $\gamma$ and $L$ are constants,
and $T \leq n$.  Applying the Hoeffding-Azuma inequality with $\lambda
= n^{-1/3}$, we find that for each $t \leq T \leq n$,
\begin{displaymath}
  \pr{X_t \geq 0}
  \ \leq \ \exp\left\{-\frac{n^{-2/3}}{2 \cdot \left( \frac{2}{n} (C \log n)^2 \right)^2 t}\right\}
  \ \leq \ \exp\left\{-\frac{n^{1/3}}{8 C^4 \log^4 n}\right\}.
\end{displaymath}
A union bound over all $t \leq T$ shows that \whp, all $X_t
< 0$.  Furthermore, Fact \ref{fact:SW} implies that \whp, $G_T$ has a
$K', c'$ component tail.

To complete our argument, we claim that whenever all of these
high-probability events happen, then all $\mathcal{E}_t$ occur for $0 \leq t
\leq T$.  We prove this by induction on $t$.  Each $\mathcal{E}_t$ has two
parts: a component tail and an upper bound on $S(G_t)$.  We will never
need to worry about checking the $K', c'$ component tail property
because $G_t \subset G_T$ for all these $t$, and we already conditioned on $G_T$
having that (monotone) property.  We concentrate on the upper bounds for
$S(G_t)$ in the remainder of this proof.  For the base case $t=0$, the
susceptibility part of $\mathcal{E}_0$ is immediate by definition since $S(G_0)
= \phi(0) < \phi(0) + e^{\frac{5L}{\gamma} \frac{0}{n}}
n^{-\frac{1}{3}}$.  For our induction step, given that $\mathcal{E}_{t-1}$
occurs, the definition of $X_t$ is then $S(G_t) - \phi(t) -
e^{\frac{5L}{\gamma} \frac{t}{n}} n^{-\frac{1}{3}}$ instead of the
alternative $X_{t-1}$.  Yet we assumed that $X_t < 0$, so that gives
the susceptibility part of $\mathcal{E}_t$, and completes the induction.

Therefore, we conclude that $\mathcal{E}_T$ occurs \whp, which in particular
means that $S(G_T) \leq \phi(T) + e^{\frac{5L}{\gamma} \frac{T}{n}}
n^{-\frac{1}{3}} = \frac{L}{\gamma} + o(1)$, as desired.  \hfill
$\Box$

\section{Offline avoidance of giants}
\label{sec:offline-avoid}

In this section, we prove Theorem \ref{thm:offline-avoid}, which has
two parts, a lower and an upper bound.  The lower bound relies on the
following relationship between orientability and decomposition.
Recall that we call a graph $r$-orientable if it is possible to orient
all edges such that all in-degrees are at most $r$.

\begin{lemma}
  \label{lem:orientable-color}
  The edges of any $r$-orientable graph $G$ can be colored with $r$
  colors such that for every pair of distinct vertices $u,v$, there
  are at most 2 monochromatic paths in each color connecting $u$ and
  $v$.
\end{lemma}

\noindent \textbf{Proof.}\, Fix an orientation of $G$ with all
in-degrees at most $r$, and greedily color the edges by $r$ colors
so that at each vertex, all incoming edges are differently colored.
Consider a particular color class.  By construction, it is a
directed graph with all in-degrees at most 1, so it is a disjoint
union of unicyclic components.  Then, every pair of vertices is
linked by at most two paths in that color, as desired.  \hfill
$\Box$

\vspace{3mm}

The previous lemma produces a coloring whose connectivity is very
fragile.  Our next lemma quantifies this, showing that the (\emph{a
  priori}, possibly large) monochromatic components shatter easily.

\begin{lemma}
  \label{lem:components-break}
  For any $\epsilon > 0$, there is $c > 0$ such that the following
  holds.  Let $G$ be a graph on $n$ vertices with maximum degree $\log n$, where
  every pair of distinct vertices is connected by at most 2 distinct
  paths.  Independently delete each edge of $G$ with probability
  $\epsilon$.  Then, \whp\ all connected components of the resulting
  graph have order at most $n e^{-c \frac{\log n}{\log \log n}} =
  o(n)$.
\end{lemma}

\noindent \textbf{Proof.}\, Define $c$ such that $1-\epsilon = e^{-8
  c}$, and let the \emph{susceptibility}\/ be $\frac{1}{n} \sum_v
C_v$, where $C_v$ is the size of the connected component containing
$v$.  Let the random variable $S$ be the susceptibility of the graph
$G'$ which remains after the edge deletions.  Since $\sum_v C_v$
equals the sum of the squares of the component sizes, all components
of $G'$ have order at most $\sqrt{nS}$.  Thus, it suffices to show
that $S \leq n e^{-2c \frac{\log n}{\log \log n}}$ \whp.

Fix an arbitrary vertex $v$.  Since $G$ has maximum degree $\log n$,
the total number of vertices within distance $D = \frac{1}{2}
\frac{\log n}{\log \log n}$ of $v$ is at most $(\log n)^D =
\sqrt{n}$.  Any other vertex $u$ has probability at most $2
(1-\epsilon)^D = 2 e^{-8 c D}$ of being connected to $v$
after the deletion.  This is because there are at most 2 paths between
$u$ and $v$, and each path has length at least $D$.  Therefore, by
linearity of expectation, the expected size of the component
containing $v$ is $\E{C_v} \leq \sqrt{n} + n \cdot 2 e^{-4c
  \frac{\log n}{\log \log n}} \leq n e^{-3c \frac{\log
    n}{\log \log n}}$.  Another application of linearity of
expectation gives $\E{S} \leq n e^{-3c \frac{\log n}{\log
    \log n}}$.  So, by Markov's inequality, $S$ exceeds $n
e^{-2c \frac{\log n}{\log \log n}}$ with probability at most
$e^{-c \frac{\log n}{\log \log n}} = o(1)$, completing the proof.
\hfill $\Box$

\vspace{3mm}

\noindent \textbf{Remark.}\, The self-contained argument above only
requires a relatively weak maximum degree condition, and is sufficient
for our purposes.  It is worth mentioning that under the stronger
assumption that $G$ is a random graph, one can use the substantially
less trivial Lemma 11 of \cite{BK} to sharpen the eventual bound to
$n^{1-c_\epsilon}$, as Sp\"ohel, Steger and Thomas do in \cite{SST}.
Indeed, that lemma claims that if $m = cn$, then there is a large
constant $K$ such that in $G_{n,m}$, the number of vertices within
distance $(\log n)/K$ of any vertex $v$ is at most $n^{\frac{\log
    2K}{K}}$ \whp.  Using this fact above instead of our exploration
to depth $\frac{1}{2} \frac{\log n}{\log \log n}$ bounds all connected
components below $n^{1-c_\epsilon}$.

\vspace{3mm}

Let us now state the result of Cain, Sanders and Wormald
\cite{CSW}, and Fernholz and Ramachandran \cite{FR}, on the matching
thresholds for orientability and average degree.

\begin{fact}
  \label{fact:orientable-threshold}
  For any integer $r \geq 2$, there is an explicit threshold $\psi_r$
  such that the following holds.  For any $\epsilon > 0$, if $m <
  (\psi_r - \epsilon)n$, then $G_{n,m}$ is $r$-orientable \whp.  On
  the other hand, if $m > (\psi_r + \epsilon)n$, then \whp\ $G_{n,m}$
  contains a subgraph with average degree at least $2r + c_\epsilon$,
  where $c_\epsilon > 0$.
\end{fact}

We are now ready to prove the lower bound, which we first translate to
$G_{n,p}$ for convenience.  By applying Fact \ref{fact:coupling} and
rescaling $\epsilon$, it suffices to show that if $p =
2(1-\epsilon)(\psi_r - \epsilon)/n$, then \whp\ there is a coloring of
$G_{n,p}$ where every color class has all components of order $o(n)$.

\vspace{3mm}

\noindent \textbf{Proof of lower bound of Theorem
  \ref{thm:offline-avoid}.}\, Let $p' = 2(\psi_r - \epsilon)/n$, and
observe that $G_{n,p}$ can be obtained from $G' = G_{n,p'}$ by
independently deleting each edge with probability $\epsilon$.  First,
consider the graph $G'$ before deletions.  By Fact
\ref{fact:orientable-threshold}, $G'$ is $r$-orientable \whp.  Also,
it is easy to see that since $np$ is at most the constant $2 \psi_r$,
$G'$ has maximum degree at most $\log n$ \whp.  Indeed, each
individual degree is distributed as $\bin(n-1, p)$, and $\pr{\bin(n-1,
  p) > \log n} \leq \binom{n}{\log n} p^{\log n} \leq \big(
\frac{enp}{\log n} \big)^{\log n}$.  Since $np$ is bounded by a
constant, this is $o(n^{-1})$, so a union bound over all $n$ vertices
implies that the maximum degree is at most $\log n$ \whp.

Thus, by Lemma
\ref{lem:orientable-color}, we can color the edges of $G'$ so that
every pair of distinct vertices is connected by at most two paths in
each color.  This, together with our degree bound and Lemma
\ref{lem:components-break}, shows that after deleting each edge of
$G'$ independently with probability $\epsilon$ (to obtain
$G_{n,p}$), \whp\ all color classes have connected components of
order only $o(n)$. \hfill $\Box$

\vspace{3mm}

For the upper bound, we use the second half of Fact
\ref{fact:orientable-threshold}, which gives a subgraph of high
average degree.  It turns out that this is already enough to ensure a
giant.  To see this, we first show that small sets of vertices
typically induce low average degree in the random graph.

\begin{lemma}
  \label{lem:small-subsets-small-avg-deg}
  For any $\lambda, \epsilon > 0$, there is a constant $c > 0$ such
  that in $G_{n,p}$ with $p = \frac{\lambda}{n}$, \whp\ every set of
  at most $cn$ vertices induces a subgraph with average degree less
  than $2+\epsilon$.
\end{lemma}

\noindent \textbf{Proof.}\, Without loss of generality, assume that
$\epsilon < 1$ and $\lambda \geq 1$.  Let $c = (e^3
\lambda^2)^{-\frac{2}{\epsilon}}$.  We will take a union bound over
all subsets of $t \leq cn$ vertices.  For a fixed value of $t$, the
probability that some $t$-set of vertices induces at least $\big(1 +
\frac{\epsilon}{2}\big)t$ edges is at most
\begin{eqnarray*}
  \binom{n}{t} \cdot \pr{\bin\left[ \binom{t}{2}, \frac{\lambda}{n} \right] \geq \left(1 + \frac{\epsilon}{2}\right)t}
  &\leq& {n \choose t}
  \cdot {t^2/2 \choose (1 + \frac{\epsilon}{2})t} \left(\frac{\lambda}{n}\right)^{(1 + \frac{\epsilon}{2})t} \\
  &\leq& \left( \frac{en}{t} \right)^{t}
  \cdot \left( \frac{et^2/2}{(1 + \frac{\epsilon}{2})t} \cdot \frac{\lambda}{n}\right)^{(1 + \frac{\epsilon}{2})t} \\
  &=& \left[
    \left( \frac{en}{t} \right)
    \cdot
    \left( \frac{e \lambda}{2 + \epsilon} \cdot \frac{t}{n}\right)^{1 + \frac{\epsilon}{2}}
  \right]^t \\
  &=& \left[
    e \left( \frac{e \lambda}{2 + \epsilon} \right)^{1 + \frac{\epsilon}{2}}
    \cdot
    \left(
      \frac{t}{n}
    \right)^{\frac{\epsilon}{2}}
  \right]^t \\
  &\leq& \left[
    \frac{e^3 \lambda^2}{2}
    \cdot
    \left(
      \frac{t}{n}
    \right)^{\frac{\epsilon}{2}}
  \right]^t.
\end{eqnarray*}
To complete our union bound, we sum the final expression over the
range $1 \leq t \leq cn$.  We split this into two intervals,
separating at $t = \log n$.  Observe that the quantity in the square
brackets increases in $t$, and reaches $\frac{1}{2}$ when $t = cn$.
So, the sum over the interval $\log n \leq t \leq cn$ is at most
$\sum_{\log n}^{cn} 2^{-t} = o(1)$.  For the other interval $t < \log
n$, the square bracket is still at most $\frac{1}{2} \leq 1$, so we
can ignore the outer exponentiation and conclude that the final
expression is at most $\frac{e^3 \lambda^2}{2} \cdot \big( \frac{\log
  n}{n} \big)^{\frac{\epsilon}{2}}$.  Multiplying this by the number
of values of $t$ in this interval ($\log n$), we see that the final
sum is still $o(1)$.  Therefore, the property holds \whp, as claimed.
\hfill $\Box$

\vspace{3mm}

From this, we immediately derive the following useful corollary, which
ensures a giant in any subgraph of average degree at least $2 +
\epsilon$.

\begin{corollary}
  \label{cor:avg-deg-2-giant}
  For any $\lambda, \epsilon > 0$, there is a constant $c > 0$ such
  that in $G_{n,m}$ with $m = \lambda n$, \whp\ every subgraph with
  average degree at least $2+\epsilon$ contains a connected component
  of order at least $cn$.
\end{corollary}

\noindent \textbf{Proof.}\, By the previous lemma and Fact
\ref{fact:coupling}, \whp\ $G_{n,m}$ has the property that every set
of at most $cn$ vertices induces a subgraph with average degree less
than $2+\epsilon$.  Then, consider any subgraph $H$ with average
degree at least $2+\epsilon$.  Separating $H$ into its connected
components, we find that some component must have average degree at
least $2+\epsilon$.  Therefore, that component must have order at
least $cn$, as desired. \hfill $\Box$

\vspace{3mm}

\noindent \textbf{Proof of upper bound of Theorem
  \ref{thm:offline-avoid}.}\, By Fact \ref{fact:orientable-threshold},
if $m > (\psi_r + \epsilon)n$, \whp\ $G_{n,m}$ contains a subgraph $H$
with average degree at least $2r + c_\epsilon$.  No matter which
colors appear on the edges of $H$, some color class will have average
degree at least $2+c_\epsilon/r$, and therefore contain a giant \whp\
by Corollary \ref{cor:avg-deg-2-giant}.  \hfill $\Box$

\section{Online avoidance of giants}
\label{sec:online-avoid}

In this section, we consider the online case of the avoidance
problem.  We first show that a natural adaptation of the offline
algorithm gives an asymptotically sharp result for large numbers of
colors.  Then, we consider the other extreme with 2 colors, and show
that the trivial bound is not tight there.

\subsection{Many colors}

Our offline algorithm avoided giant components by orienting edges to
minimize in-degrees.  By replacing the offline orientation procedure
with an online one, this strategy naturally extends to the online
setting.  Online edge orientation has been extensively studied, in
the famous equivalent formulation known as the ``power of two random
choices'' with balls and bins (see \cite{MRS} for a survey of
results).  Indeed, that setting had $n$ bins, with $kn$ balls coming
sequentially, each with two independent random choices for a
destination bin.  The objective was to control the maximum load
across all of the bins.  This can be interpreted as a graph
orientation problem, where each pair of bin choices corresponds to
an incoming edge with the two choices as endpoints.  The edge's
orientation records which bin the ball is sent to, and the goal of
controlling the maximum in-degree is precisely the same as that of
controlling the maximum load in the balls-and-bins problem.

It is now well-known that when the objective is to minimize the
maximum in-degree, the stochastically optimal online orientation
strategy is to always orient each incoming edge towards the endpoint
which currently has lower in-degree.  However, it turns out that for
the purpose of proving Theorem \ref{thm:online-avoid}, one can use a
random orientation strategy, which is easier to analyze.  Our
coloring algorithm, which we call \orient, internally maintains a
set of orientations for all edges it has seen.  To color a new edge
$e$, it randomly orients it with equal probability toward one of its
endpoints.  Let the new in-degree of that endpoint be $d$.  If $d <
r$, then color $d$ is used for the edge $e$.  Otherwise, color $r$
is used.  Observe that just as in Lemma \ref{lem:orientable-color},
each of the first $r-1$ color classes is a disjoint union of
unicyclic components.  Therefore, each of these color classes has
every pair of vertices connected by at most two paths, so it will
shatter by the same argument as in the proof of Theorem
\ref{thm:offline-avoid}.

The new challenge in this section is to control the $r$-th color
class.  Fortunately, it turns out that it is extremely sparse.  To
prove this, it is more convenient to work in the random directed
graph $\Gnpd$, in which each of the $n(n-1)$ possible directed edges
appears independently with probability $p/2$.  Note that in this
model, it is possible for both $\overrightarrow{uv}$ and
$\overleftarrow{uv}$ to be present simultaneously.  Our first claim
is that $\Gnpd$ typically has no long cycles containing many
vertices of high in-degree.  This is relevant because every edge in
color $r$ has an endpoint with in-degree at least $r$.

\begin{lemma}
  \label{lem:no-long-cycles}
  For any $\epsilon > 0$, the following holds for every sufficiently
  large constant $r$.  Let $\Gd = \Gnpd$ be a random directed
  graph with $p = (1-\epsilon)\frac{2r}{n}$, and $G$ be the undirected
  graph on the same vertex set obtained by collapsing all edges
  between each vertex pair into a single undirected edge.  Then, \whp\
  $G$ does not contain any cycles of length at least $\sqrt[4]{\log
    n}$ for which at least half of the vertices on the cycle had
  in-degree at least $r$ in $\Gd$.
\end{lemma}

\noindent \textbf{Proof.}\, We will use a union bound to show that a
large family of objects do not appear in the random directed graph.
Let us define an \emph{isomorphism type}\/ to be a cycle, say with
vertices $v_1, \ldots, v_t$, whose edges have been oriented
arbitrarily, along with a subset of at least $t/2$ of its vertices
which have been designated as ``high-in-degree vertices.''  Note that
we do not require the edges of the cycle to be oriented in a
consistent direction.  The number of distinct $t$-vertex isomorphism
types is at most $2^t \cdot 2^t$, because each of the $t$ edges can be
oriented in 2 ways, and the number of different subsets of vertices
that can be designated as high-in-degree is at most $2^t$.

We say that $\Gnpd$ contains a copy of this isomorphism type if there
is an embedding of the vertices $v_i$ such that all consecutive edges
$v_i v_{i+1}$ are present in the correct direction, and all designated
high-in-degree vertices $v_i$ already have in-degree at least $r-2$
from vertices other than $v_{i-1}, v_{i+1}$.  We do not restrict our
attention to induced copies, so other edges may also be present.  If
we can show that over all isomorphism types with $t \geq \sqrt[4]{\log
  n}$, the expected total number of copies in $\Gnpd$ is $o(1)$, then
we will be done by Markov's inequality.

So, let us focus on a particular isomorphism type with $t$ vertices.
There are at most $n^t$ ways to embed the $t$ vertices of the cycle.
Each edge $v_i v_{i+1}$ independently appears with its correct
orientation with probability exactly $p/2$.  Next, consider a
designated high-in-degree vertex $v_i$.  Crucially, we only require
in-degree at least $r-2$ from vertices \emph{other than}\/ $v_{i-1}$
and $v_{i+1}$.  The reason for this exclusion is that the previous
step may already have exposed the edges $\overrightarrow{v_{i-1} v_i}$
and $\overrightarrow{v_{i+1} v_i}$.  But now, since our model even
allows edges in both directions between vertex pairs, the probability
that each designated vertex indeed has high in-degree is independently
$\pr{\bin\big[n-3, \frac{p}{2}\big] \geq r-2}$.  Since $p =
(1-\epsilon)\frac{2r}{n}$ and $r$ is large, each of these individual
probabilities is bounded by the probability that $\bin\big[n-3,
(1-\epsilon)\frac{r}{n}\big]$ exceeds its mean by at least an
$\frac{\epsilon}{2}$-fraction.  By the Chernoff bound, this happens
with probability at most $e^{-c_\epsilon r}$ for some constant
$c_\epsilon$.  By choosing large enough $r$, we may assume that this
is below $\frac{1}{64r^2}$.  Putting everything together, we find that
the expected number of copies of the isomorphism type in $\Gnpd$ is at
most
\begin{displaymath}
  n^t \left( \frac{p}{2} \right)^t \left( \frac{1}{64r^2} \right)^{t/2}
  \ \leq \
  \left( \frac{1}{8} \right)^t.
\end{displaymath}
We initially showed that the number of distinct $t$-vertex isomorphism
types is at most $4^t$, so the expected total number of copies of all
$t$-vertex isomorphism types is at most $2^{-t}$.  This is a geometric
series, so its sum over all $t \geq \sqrt[4]{\log n}$ is still $o(1)$,
as desired.  \hfill $\Box$

\vspace{3mm}

\noindent \textbf{Remark 1.}\, Since we had a convergent geometric
series at the end of the proof, the $\sqrt[4]{\log n}$ bound is not
tight.  In fact, any function which grows with $n$ is sufficient.

\vspace{3mm}

\noindent \textbf{Remark 2.}\, If one is interested in beating the
trivial bound, which corresponds to $p \approx \frac{r}{n}$, one can
choose $\epsilon$ to be extremely close to, but just below,
$\frac{1}{2}$.  One can numerically check that if $\epsilon = 0.4999$
and $r \geq 51$, then the probability that $\bin\big[n,
(1-\epsilon)\frac{r}{n}\big]$ exceeds $r-2$ is at most $\frac{1}{16.1
  r^2}$ for large $n$, because the Binomial converges to a Poisson
variable with mean $0.5001 r$.  Continuing the argument, this will
show that the expected number of appearances of all $t$-vertex
isomorphism types is at most $\big( \frac{4}{\sqrt{16.1}} \big)^t$,
which is still a convergent geometric series, so the same result will
follow.

\vspace{3mm}

Next, we establish an easy bound which holds for ordinary random
graphs.

\begin{lemma}
  \label{lem:small-components-not-complex}
  For every constant $c$, \whp\ in $G_{n,p}$ with $p = \frac{c}{n}$, every
  set of $t \leq \sqrt[3]{\log n}$ vertices induces at most $t$ edges.
\end{lemma}

\noindent \textbf{Proof.}\,
The expected number of sets with $t\leq \sqrt[3]{\log n}$ and at least $t+1$ edges
can be bounded by
$$\sum_{t=4}^{\sqrt[3]{\log n}}\binom{n}{t}\binom{\binom{t}{2}}{t+1}p^{t+1}\leq \sum_{t=4}^{\sqrt[3]{\log n}}
\frac{t}{en}\left(\frac{ne}{t}\cdot \frac{tec}{2n}\right)^{t+1}=o(1).$$

\hfill $\Box$

\vspace{3mm}

We now combine the previous two lemmas to show that the $r$-th color
class shatters easily.  In the proof of Lemma
\ref{lem:components-break}, the control of connectivity was done by
bounding the number of distinct paths between every pair of
vertices. This time, we use the notion of an \emph{essential edge}.
We say that an edge $e$ on a path is \emph{essential} if every other
path connecting the same endpoints also contains $e$.  It turns out
that in the $r$-th color class, every long path contains a huge
number of essential edges.

\begin{lemma}
  \label{lem:r-fragile}
  For any $\epsilon > 0$, the following holds \whp\ for every
  sufficiently large constant $r$.  Let $G$ be the graph formed by the
  $r$-th color class after $(1-\epsilon)rn$ independent random edges
  have been colored by \orient.  Then every path in $G$ of length at
  least $\sqrt[3]{\log n}$ has the property that more than half of its
  edges are essential.
\end{lemma}

\noindent \textbf{Proof.}\, Since each (random) incoming edge is
randomly directed by \orient, one can think of the input as a sequence
of random directed edges, which is then deterministically colored
using the rule in \orient.  By a similar argument to Fact
\ref{fact:coupling}, it suffices to consider the more convenient model
where the input sequence is a random permutation of the edges of a
random directed graph $G = \Gnpd$ with $p =
(1-\epsilon)\frac{2r}{n}$.

Note that if an edge of $G$ is oriented toward a vertex with in-degree
less than $r$, then regardless of the permutation, it will never be
colored $r$.  So, let $H \subset G$ be obtained by deleting all edges
oriented into vertices of in-degree less than $r$.  Then $H$ entirely
contains the $r$-th color class.  Let $A$ be the set of vertices whose
in-degrees were less than $r$, and let $B$ be those that had in-degree
at least $r$.  Observe that we deleted all edges oriented toward
vertices in $A$, so $A$ spans no edges in $H$.  In particular, any
cycle in $H$ has at least half of its vertices in $B$, i.e., with
in-degree at least $r$.

Therefore, by Lemma \ref{lem:no-long-cycles}, \whp\ all cycles in $H$
have length at most $\sqrt[4]{\log n}$.  Also, condition on the result
of Lemma \ref{lem:small-components-not-complex}, which shows that in
$G$ (and hence also $H$), every set of $t \leq \sqrt[3]{\log n}$
vertices induces at most $t$ edges.  These two graph properties will
be enough to show that long paths in $H$ contain many essential edges.

Let $P = v_1, \ldots, v_t$ be a path in $H$ with length at least
$\sqrt[3]{\log n}$.  Suppose for contradiction that at least half of
its edges are non-essential.  We claim that since $\sqrt[4]{\log n}
\ll \sqrt[3]{\log n}$, there must be non-essential edges $v_i v_{i+1}$
and $v_j v_{j+1}$ such that $i < j$ and $3 \sqrt[4]{\log n} <
j-i < 7 \sqrt[4]{\log n}$.  Indeed, if this were false, then out of
the $7 \sqrt[4]{\log n}$ edges immediately following each
non-essential edge in $P$, at least $\frac{4}{7}$-fraction of them
would be essential.  Then an averaging argument would contradict the
fact that at least half of the edges were non-essential.

Now, since $v_i v_{i+1}$ is non-essential, there is another path $P' =
w_1, \ldots, w_s$ with $w_1 = v_1$ and $w_s = v_t$ which avoids the
edge $v_i v_{i+1}$.  Let $a$ be the largest index such that $w_a \in
\{v_1, \ldots, v_i\}$, and let $b$ be the next index after $a$ such
that $w_b \in P$.  These exist because $P$ and $P'$ both contain $v_1$
and $v_t$.  Note that by definition, $w_b$ is actually in $\{v_{i+1},
\ldots, v_t\}$, and the segment of $P'$ from $w_a$ to $w_b$ intersects
$P$ only at $w_a$ and $w_b$.  So, there is a cycle $C_1$ formed by
going from $w_a$ to $w_b$ along $P'$, and then back to $w_a$ along
$P$.  Importantly, the common edges between $C_1$ and $P$ are a
contiguous interval containing the edge $v_i v_{i+1}$.

Similarly, we can find a cycle $C_2$ containing the edge $v_j
v_{j+1}$.  Crucially, $C_1$ and $C_2$ are distinct (although not
necessarily disjoint) because $j-i > 3 \sqrt[4]{\log n}$ and we
conditioned on all cycles being shorter than $\sqrt[4]{\log n}$. Yet
$j-i < 7 \sqrt[4]{\log n}$, so the union of $C_1$, $C_2$, and the
path $v_i v_{i+1} \ldots v_j$ forms a subgraph of order $k \leq 9
\sqrt[4]{\log n}$ which spans at least $k+1$ edges. Since we also
conditioned on all such subgraphs having order at least
$\sqrt[3]{\log n}$, this is a contradiction.  Therefore, the path $P$ must have
had at least half of its edges essential.  \hfill $\Box$

\vspace{3mm}

The previous lemma shows that the $r$-th color class is typically
quite fragile as well.  We now combine this with an adaptation of
our offline argument, and prove that it is possible to avoid giants
in all colors for nearly $rn$ rounds \whp.

\vspace{3mm}

\noindent \textbf{Proof of Theorem \ref{thm:online-avoid}.}\, By
rescaling $\epsilon$, it suffices to give a randomized coloring
algorithm that avoids giants in all colors \whp, for a sequence of $m
= (1 - \epsilon)^3 rn$ independent random edges (possibly with
repetitions).  As in our proof of Theorem \ref{thm:offline-avoid}, it
is convenient to color a slightly denser random graph, because the
deletion of fictitious edges shatters all large components.

Strictly speaking, we cannot simply apply \orient\ to a larger
sequence of edges, because for this problem the input is a sequence
of $m$ edges, which must be processed online.  We will therefore
take some care in specifying how we randomly interleave the input
into a longer sequence of edges, so that all operations are clearly
online. Let us denote the final sequence of real and fictitious
edges by $e_1, \ldots, e_{m'}$, where $m' = (1-\epsilon)^2rn$.
Initially, we select a random subset of $m$ of the $m'$ indices to
correspond to the positions of the real edges.  We then generate
independent random edges for all other $e_i$, and pass the resulting
sequence to \orient. Note that since the input distribution is
uniform over all sequences of $m$ edges, the augmented sequence of
edges consists of $m'$ independent random edges.

Let $\sigma$ denote the colored sequence of $m'$ edges produced by
\orient.  The graph formed by $\sigma$ has maximum degree at most
$\log n$ \whp\ by the same argument as in the offline case.  We also
know by construction that there are at most 2 paths between every
pair of vertices in each of the first $r-1$ color classes.  For the
$r$-th color class, Lemma \ref{lem:r-fragile} ensures that \whp, all
paths longer than $\sqrt[3]{\log n}$ have at least half of their
edges essential.  Let $\mathcal{P}$ denote the collection of all of
these properties.  We will write $\sigma \in \mathcal{P}$ when all
of them hold.

Now, we delete the $\epsilon m'$ fictitious edges to recover the
coloring of the original edges.  Note that since the algorithm knows
which $m$ edges are real (that was the input), the edges to delete are
completely determined.  But crucially, it used an independent source
of randomness to interleave the original $m$ edges into the full
sequence of $m'$ edges.  Therefore, if we only condition on $\sigma$
(and not on the input), then the distribution of which $m$ edges were
original is uniform over all possible subsets of $m$ positions.
Formally, we are calculating the probability of success by summing
over all colored sequences $\sigma$ of $m'$ edges.  We have
\begin{displaymath}
  \pr{\text{success}}
  \ = \
  \sum_{\sigma} \pr{\text{success} \mid \sigma} \pr{\sigma}
  \ \geq \
  \sum_{\sigma \in \mathcal{P}} \pr{\text{success} \mid \sigma} \pr{\sigma}
\end{displaymath}
Since we showed that $\sigma \in \mathcal{P}$ \whp, it suffices to
show that $\pr{\text{success} \mid \sigma} \geq 1-o(1)$ for all
$\sigma \in \mathcal{P}$.  We noted above that conditioned on
$\sigma$, the $\epsilon m'$ edges to delete were uniformly
distributed over all subsets.  Therefore, it remains to show that
given any coloring with property $\mathcal{P}$, the deletion of a
random $\epsilon$-fraction of its edges \whp\ shatters all large
connected components.  We accomplish this by deleting every edge
independently with probability $\frac{\epsilon}{2}$, which will
imply the result by a similar coupling argument to Fact
\ref{fact:coupling}, since $\bin \big(m', \frac{\epsilon}{2} \big)
\leq \epsilon m'$ \whp.

For each of the first $r-1$ color classes, Lemma
\ref{lem:components-break} shows that all components shatter to
$o(n)$ \whp, as in the offline proof.  For the $r$-th color class,
we now adapt the proof of Lemma \ref{lem:components-break} to use
essential edges.  Indeed, let us bound the expected size of the
component $C_v$ containing a particular vertex $v$ after the
deletions.  Since the maximum degree in $G_{n,rn}$ is $\log n$, the
total number of vertices within distance $D = \frac{1}{2} \frac{\log
n}{\log \log n}$ of $v$ is at most $(\log n)^D = \sqrt{n}$.  Any
other vertex $u$ is at distance at least $D \gg \sqrt[3]{\log n}$
away from $v$, so a shortest path from $v$ to $u$ contains at least
$D/2$ essential edges.  The deletion of any essential edge
disconnects $u,v$, so if edges are deleted with probability
$\frac{\epsilon}{2}$, then $u$ and $v$ remain connected only with
probability at most $(1-\frac{\epsilon}{2})^{D/2} = e^{-c
  \frac{\log n}{\log \log n}}$ for some constant $c$.  Hence the
expected size of $C_v$ is at most $\sqrt{n} + n e^{-c \frac{\log
    n}{\log \log n}} = o(n)$, and by linearity of expectation, the
expected susceptibility $\E{S}$ of the graph after deletions is
$o(n)$.  Since the size of the largest component is at most
$\sqrt{nS}$, Markov's inequality implies that the $r$-th color class
also has all components smaller than $o(n)$, completing the proof.
\hfill $\Box$

%
%
%

\subsection{Two colors}

The trivial algorithm, which randomly colors each edge blue or red,
clearly lasts for $(1-\epsilon) n$ rounds \whp.  We now present a
better algorithm which lasts for $1.06n$ rounds \whp.  To color a new
edge $e$, it considers the set of colors $C$ that appear on isolated
edges which are incident with any of its endpoints.  (If $e$ is not
incident to any isolated edges, then $C$ is empty.)  When $C$ contains
exactly one color, the algorithm colors the edge $e$ with the other
color.  Otherwise, it randomly colors $e$ either blue or red with
equal probability.

We analyze this by tracking a certain partition of the vertex set.
Split the set of isolated edges into two groups based on their
color, and call them the \emph{red matching}\/ and the \emph{blue
matching}, respectively.  After the $k$-th round, let:
\begin{eqnarray*}
  I_k &=& \text{number of isolated vertices},\\
  B_k &=& \text{number of vertices in the blue matching}, \\
  R_k &=& \text{number of vertices in the red matching},
\end{eqnarray*}
and let $J_k = n - I_k - B_k - R_k$ be the number of remaining
vertices.  These parameters correspond to the decomposition of the
graph into its isolated vertices, the blue matching, the red matching,
and the remainder.

\begin{lemma}
  \label{lem:track-partition}
  With probability $1-o(1)$, the following hold for all $t \leq 1.1$:
  \begin{eqnarray*}
    \textstyle \big|\frac{1}{n} I_{tn} - e^{-2t}\big| &\leq& n^{-1/3}, \\
    \textstyle \big|\frac{1}{n} B_{tn} - t e^{-4t}\big| &\leq& e^8 n^{-1/3}, \\
    \textstyle \big|\frac{1}{n} R_{tn} - t e^{-4t}\big| &\leq& e^8 n^{-1/3}. \\
  \end{eqnarray*}
\end{lemma}

\noindent \textbf{Proof.}\, The probability that a particular vertex
is not incident to any of the first $tn$ edges is exactly $\big(
\frac{n-1}{n} \cdot \frac{n-2}{n-1}\big)^{tn} = \big(1 -
\frac{2}{n}\big)^{tn}$, which tends to $e^{-2t}$ from below as $n$
grows.  Routine calculus easily bounds the convergence rate by
$O(n^{-1})$, so $\E{\frac{1}{n}I_{tn}} = e^{-2t} + O(n^{-1})$.
Now consider the edge-exposure martingale where $Y_k$ is the
conditional expectation of $I_{tn}$ given the first $k$ rounds.
Changing the outcome of any particular round can only affect $I_{tn}$
by at most 2, and there are $tn$ rounds to determine $I_{tn}$, so by
the Hoeffding-Azuma inequality (see Theorem 7.4.1 of \cite{AS})
$I_{tn}$ is within (say) $\frac{1}{2}n^{2/3}$ of its expectation with
probability $e^{-\Omega(n^{1/3})}$.  This gives the desired asymptotic
for $I_{tn}$.

We estimate $B_{tn}$ next.  We claim that conditioned on the first $k$
incoming edges $e_1, \ldots, e_k$, the expected change $B_{k+1} - B_k$
is
\begin{equation}
  \label{eq:deltaB}
  \E{B_{k+1}-B_k \mid e_1, \ldots, e_k}
  \ = \
  2 \cdot \left( \frac{I_k}{n} \right)^2 \frac{1}{2}
  - \frac{4B_k}{n}
  +
  O(n^{-1}).
\end{equation}
The first summand comes from the creation of a blue isolated edge from
2 isolated vertices, which contributes 2 to $B_k$.  The probability
that both endpoints are isolated vertices is $\frac{I_k}{n} \cdot
\frac{I_k - 1}{n-1}$.  Since $\frac{1}{n(n-1)} - \frac{1}{n^2} =
O(n^{-3})$ and $I_k \leq n$, this is $\big(\frac{I_k}{n}\big)^2 -
O(n^{-1})$.  The $\frac{1}{2}$ factor comes from the fact that the
edge is randomly colored blue or red.

For the second summand, the only way we lose blue isolated edges is
when an endpoint of the incoming edge is incident to a blue isolated
edge.  The probability that the two endpoints hit two different blue
isolated edges (hence contributing $-4$) is $\frac{B_k}{n} \cdot
\frac{B_k - 2}{n-1}$.  On the other hand, the probability that they
hit exactly one isolated edge (hence contributing $-2$) is $2 \cdot
\frac{B_k}{n} \big(1 - \frac{B_k-1}{n-1} \big)$.  Thus the expected
contribution from these losses is
\begin{displaymath}
  (-4) \cdot \frac{B_k}{n} \cdot
  \frac{B_k - 2}{n-1}
  + (-2) \cdot 2 \cdot \frac{B_k}{n} \left(1 - \frac{B_k-1}{n-1} \right)
  \ = \
  -\frac{4B_k}{n} + O(n^{-1}),
\end{displaymath}
matching the second summand.

Since we showed that $\frac{1}{n} I_{tn} = (1-o(1))e^{-2t}$ \whp,
equation \eqref{eq:deltaB} suggests that $b(t)=\frac{1}{n}B_{tn}$
should satisfy the differential equation
\begin{displaymath}
  \frac{db}{dt} \ = \ (e^{-2t})^2 - 4b, \quad \quad b(0) \ = \ 0,
\end{displaymath}
whose solution is $b(t) = t e^{-4t}$.

We now verify this formally, using the same method as for the proof of
Theorem \ref{thm:evolve-suscep}.  For each $k$, let $\mathcal{E}_k$ be
the event that $\big|\frac{1}{n} I_k -e^{-\frac{2k}{n}}\big| \leq
n^{-\frac{1}{3}}$ and $\big|\frac{1}{n} B_k -b\big(\frac{k}{n}\big)\big|
\leq e^{\frac{7k}{n}}n^{-\frac{1}{3}}$.  Now, consider the sequence of random
variables
\begin{displaymath}
  W_k \ = \ \begin{cases}
    B_k-nb\big(\frac{k}{n}\big)-e^{\frac{7k}{n}}n^{\frac{2}{3}} & \text{if }\cE_{k-1}\text{ occurs}, \\
    W_{k-1} & \text{otherwise}.
  \end{cases}
\end{displaymath}
We claim that $W_k$ is a supermartingale.
Assume that $\cE_k$ occurs. Then, using \eqref{eq:deltaB} we obtain
\begin{align*}
&\E{W_{k+1}-W_k \mid e_1, \ldots, e_k, \mathcal{E}_k}\\
&\alignsp \leq \ \left( \frac{I_k}{n} \right)^2  - \frac{4B_k}{n}  +  O(n^{-1})
-n\left[b\left(\frac{k+1}{n}\right)-b\left(\frac{k}{n}\right)\right]
-\left[e^{\frac{7(k+1)}{n}}-e^{\frac{7k}{n}}\right]n^{2/3}.
\end{align*}
Since $\frac{I_k}{n} \leq e^{-\frac{2k}{n}} + n^{-\frac{1}{3}}$ and
$e^{-\frac{2k}{n}} \leq 1$, we have $\big(\frac{I_k}{n}\big)^2 \leq
e^{-\frac{4k}{n}} + 2 n^{-\frac{1}{3}} + O(n^{-\frac{2}{3}})$.
Similarly, $-\frac{4B_k}{n} \leq -4b\big(\frac{k}{n}\big) + 4
e^{\frac{7k}{n}} n^{-\frac{1}{3}}$.  Recall that for any
twice-differentiable function $f$, Taylor's formula ensures that for
any $t,h$, there is some $0 \leq \xi \leq 1$ such that $f(t+h) - f(t)
= f'(t)h + \frac{1}{2} f''(t+\xi h) h^2$.  Since the second derivative
of our function $b(t)$ is bounded on the interval $0 \leq t \leq 1.1$,
Taylor's formula gives $b\big(\frac{k+1}{n}\big) -
b\big(\frac{k}{n}\big) = \frac{1}{n} b'\big(\frac{k}{n}\big) + O(n^{-2})$.  By
a similar argument, $e^{\frac{7(k+1)}{n}}-e^{\frac{7k}{n}} =
\frac{7}{n} e^{\frac{7k}{n}} + O(n^{-2})$.  Combining all of these
estimates and using $b' = e^{-4t} - 4b$, we obtain
\begin{align*}
&\E{W_{k+1}-W_k \mid e_1, \ldots, e_k, \mathcal{E}_k}\\
&\alignsp \leq \ e^{-\frac{4k}{n}} + 2 n^{-\frac{1}{3}} + O(n^{-\frac{2}{3}})
- 4b\left(\frac{k}{n}\right) + 4 e^{\frac{7k}{n}} n^{-\frac{1}{3}}
- b'\left(\frac{k}{n}\right)
- 7e^{\frac{7k}{n}} n^{-\frac{1}{3}} \\
&\alignsp = \  \left(
  2 - 3e^{\frac{7k}{n}}
  \right) n^{-\frac{1}{3}}
  + O(n^{-\frac{2}{3}}). \\
&\alignsp < \ 0,
\end{align*}
so $W_k$ is indeed a supermartingale.  Next, we bound the stepwise
differences $W_{k+1} - W_k$.  The change in $B_k$ is at most 4, and
our Taylor estimates show that the error term
$nb\big(\frac{k}{n}\big)-e^{\frac{7k}{n}}n^{\frac{2}{3}}$ changes by
at most an absolute constant because $b'\big(\frac{k}{n}\big)$ is
bounded on $k \leq 1.1n$.  Therefore, the Hoeffding-Azuma inequality
implies that since $W_0 = -n^{\frac{2}{3}}$,
\begin{equation}\label{conc1}
\pr{\exists k\leq 1.1n:W_k \geq 0}\leq e^{-\Omega(n^{1/3})}.
\end{equation}
Similarly, if
\begin{displaymath}
  \hat{W}_k \ = \ \begin{cases}
    B_k-nb\big(\frac{k}{n}\big)+e^{\frac{7k}{n}}n^{\frac{2}{3}} & \text{if }\cE_{k-1}\text{ occurs},\\
    \hat{W}_{k-1} & \text{otherwise}.
  \end{cases}
\end{displaymath}
then
\begin{align*}
&\E{\hat{W}_{k+1}-\hat{W}_k \mid e_1, \ldots, e_k, \mathcal{E}_k}\\
&\alignsp \geq \
\left( \frac{I_k}{n} \right)^2  - \frac{4B_k}{n}  +  O(n^{-1})
-n\left[b\left(\frac{k+1}{n}\right)-b\left(\frac{k}{n}\right)\right]
+\left[e^{\frac{7(k+1)}{n}}-e^{\frac{7k}{n}}\right]n^{2/3}.
\end{align*}
Since $\frac{I_k}{n} \geq e^{-\frac{2k}{n}} - n^{-\frac{1}{3}}$ and
$e^{-\frac{2k}{n}} \leq 1$, we have $\big(\frac{I_k}{n}\big)^2 \geq
e^{-\frac{4k}{n}} - 2 n^{-\frac{1}{3}}$.  Also, $-\frac{4B_k}{n} \geq
-4b\big(\frac{k}{n}\big) - 4 e^{\frac{7k}{n}} n^{-\frac{1}{3}}$.
Using the same estimates as before for $b\big(\frac{k+1}{n}\big) -
b\big(\frac{k}{n}\big)$ and $e^{\frac{7(k+1)}{n}}-e^{\frac{7k}{n}}$, we obtain
\begin{align*}
&\E{\hat{W}_{k+1}-\hat{W}_k \mid e_1, \ldots, e_k, \mathcal{E}_k}\\
&\alignsp \geq \ e^{-\frac{4k}{n}} - 2 n^{-\frac{1}{3}}
- 4b\left(\frac{k}{n}\right) - 4 e^{\frac{7k}{n}} n^{-\frac{1}{3}}
+ O(n^{-1})
- b'\left(\frac{k}{n}\right)
+ 7e^{\frac{7k}{n}} n^{-\frac{1}{3}} \\
&\alignsp = \ \left(
  -2 + 3e^{\frac{7k}{n}}
  \right) n^{-\frac{1}{3}}
  + O(n^{-1}). \\
&\alignsp > \ 0,
\end{align*}
so $\hat{W}_k$ is a submartingale.  Applying the Hoeffding-Azuma
inequality once again we see that
\begin{equation}\label{conc2}
\pr{\exists k\leq 1.1n:\hat{W}_k\leq 0}\leq e^{-\Omega(n^{1/3})}.
\end{equation}
We have now shown that \whp, $W_k < 0$, $\hat{W}_k > 0$, and
$\big|\frac{1}{n} I_k - e^{-\frac{2k}{n}}\big| \leq n^{-1/3}$ for every $k \leq
1.1n$.  Whenever these all happen, the same induction argument as in
the conclusion of the proof of Theorem \ref{thm:evolve-suscep} shows
that every $\mathcal{E}_k$ necessarily holds as well.  In particular,
\begin{displaymath}
\left| B_k - nb\left(\frac{k}{n}\right) \right|
\ \leq \
e^{\frac{7k}{n}} n^{\frac{2}{3}}
\ < \
e^8 n^{\frac{2}{3}},
\end{displaymath}
for all $k \leq 1.1n$.  This completes the proof for $B_{tn}$, and the
result for $R_{tn}$ follows by symmetry.  \hfill $\Box$

\vspace{3mm}

Now that we have control of the vertex partition, we study the
evolution of the susceptibility.  We have symmetry between blue and
red, so it suffices to show that the susceptibility of the blue color
class does not blow up before $1.06n$ rounds.  Let $X_k$ be the sum of
the squares of the component sizes in the blue color class after the
$i$-th round.  Note that this is precisely $n$ times the
susceptibility of the blue color class.  We claim that $\frac{1}{n}
X_{tn}$ tracks $x(t)$, which is the solution of the differential
equation
\begin{equation}\label{xdiff}
\frac{dx}{dt} \ = \ x^2 + 3b^2 - 2b x, \quad \quad x(0) \ = \ 1,
\end{equation}
where $b(t) = t e^{-4t}$.  Numerical methods confirm that this
differential equation blows up only at $t \approx 1.065$, and $x(t)
\leq 209$ for all $t \leq 1.06$.

\begin{lemma}
  \label{lem:deltaX}
  Suppose that $\frac{1}{n} X_k < 210$.  Then the expected change in
  $X_k$ is:
  \begin{displaymath}
    \E{X_{k+1} - X_k \mid e_1, \ldots, e_k; {\textstyle \frac{1}{n} X_k < 210}}
    \ = \
    \left( \frac{X_k}{n} \right)^2 +
    \frac{1}{n^2} \left[
      4B_k^2 - 4B_k X_k - R_k^2 + 2R_k X_k
    \right]
    + O(n^{-1}).
  \end{displaymath}
\end{lemma}

\noindent \textbf{Proof.}\, Let the connected components in the blue
color class be $C_1$, $C_2$, \ldots. Suppose that the $(k+1)$-st
edge has endpoints in $C_i, C_j$.  If $i=j$, or if the edge is colored
red, then the sum of the squares of the blue components does not
change.  Otherwise, it increases by exactly $(|C_i| + |C_j|)^2 -
|C_i|^2 - |C_j|^2 = 2 |C_i| |C_j|$.  Therefore,
\begin{displaymath}
  \E{X_{k+1} - X_k \mid e_1, \ldots, e_k; {\textstyle \frac{1}{n} X_k < 210}}
  \ = \
  \sum_{i \neq j} 2|C_i| |C_j| \cdot \frac{|C_i|}{n} \frac{|C_j|}{n-1} \cdot p_{ij}
\end{displaymath}
where $p_{ij}$ is the probability that an with endpoints
in $C_i$ and $C_j$ is colored blue.  Note that $p_{ij}$ is usually
$\frac{1}{2}$, but is sometimes 0 or 1 when the endpoints hit
isolated edges.  The factor of $n-1$ in the denominator is
cumbersome, so we will replace it with an $n$.  To do this, note
that $\sum_{i \neq j} 2|C_i|^2 |C_j|^2 \cdot p_{ij} \leq 2 ( \sum_i
|C_i|^2 )^2 = 2X_k^2 \leq 2 (210 n)^2 = O(n^2)$. Since
$\frac{1}{n(n-1)} - \frac{1}{n^2} = O(n^{-3})$, the total additive
error we will make by replacing the $n-1$ with an $n$ is
$O(n^{-1})$.  Therefore,
\begin{displaymath}
  \E{X_{k+1} - X_k \mid e_1, \ldots, e_k; {\textstyle \frac{1}{n} X_k < 210}}
  \ = \
  \frac{2}{n^2} \sum_{i \neq j} |C_i|^2 |C_j|^2 \cdot p_{ij} + O(n^{-1}).
\end{displaymath}
Let $S$ be the right hand side of the final equality, and let $S'$ be
what it would be if all $p_{ij}$ were equal to $\frac{1}{2}$.  Then
\begin{equation}\label{upsus1}
  S'
  \ = \
  \frac{1}{n^2} \sum_{i \neq j} |C_i|^2 |C_j|^2 + O(n^{-1})
  \ \leq \
  \left( \frac{X_k}{n} \right)^2 + O(n^{-1}).
\end{equation}
Now we estimate the total error we made in $S'$ by replacing all
$p_{ij}$ with $\frac{1}{2}$.  Whenever $p_{ij} = 0$, we overestimated
by $\frac{1}{n^2} |C_i|^2 |C_j|^2$, and when $p_{ij} = 1$, we
underestimated by that same amount.  To systematically examine all of
the cases when $p_{ij} \neq \frac{1}{2}$, we classify the components
$C_i$ of the blue color class into \emph{types}, which we represent
with the letters B, R, I, and J.  We say that $C_i$ has type B if it is
part of the blue matching (hence a single edge), type R if it is part
of the red matching (hence a single vertex), type I if it is an
isolated vertex, and type J otherwise.  Now we break into cases
depending on the types of $C_i$ and $C_j$.  In each case, we calculate
the sum of all $|C_i|^2 |C_j|^2$ of that type.

\begin{description}
\item[Case BB.] In this case, both $C_i$ and $C_j$ have type B,
  meaning that they are isolated edges from the blue matching.  If the
  incoming edge has one endpoint in $C_i$ and one endpoint in $C_j$,
  our algorithm will definitely color it red, so $p_{ij} = 0$.  Any
  $|C_i|^2 |C_j|^2$ of this type is precisely $2^2 \cdot 2^2 = 16$.
  The number of $C_i$ of type B is $\frac{B_k}{2}$, because the blue
  matching consists of $\frac{B_k}{2}$ isolated blue edges.  So, the
  number of pairs $C_i, C_j$ of type BB with $i \neq j$ is
  $\frac{B_k}{2} \cdot \big(\frac{B_k}{2} - 1\big) = \frac{B_k^2}{4} -
  O(n)$.  Therefore, the sum of all $|C_i|^2 |C_j|^2$ of this type is
  $4 B_k^2 - O(n)$.

\item[Cases BI, IB.] Again $p_{ij} = 0$.  Any $|C_i|^2 |C_j|^2$ of
  this type is precisely $2^2 \cdot 1^2 = 4$.  There are
  $\frac{B_k}{2} \cdot I_k$ pairs $C_i, C_j$ of type BI, and the same
  number of type IB, so the sum is $4 B_k I_k$.

\item[Cases BJ, JB.] Again $p_{ij} = 0$.  Let $Z$ be the set of
  indices $j$ such that $C_j$ has type J.  Since there are
  $\frac{B_k}{2}$ components $C_i$ of type B, the sum of $|C_i|^2
  |C_j|^2$ over all pairs of type BJ alone is
  \begin{eqnarray*}
    \frac{B_k}{2} \sum_{j \in Z} 2^2 \cdot |C_j|^2 &=& 2B_k \sum_{j \in Z} |C_j|^2 \\
    &=& 2B_k (X_k - \sum_{j \not \in Z} |C_j|^2 ) \\
    &=& 2B_k \left(X_k - I_k - R_k - \frac{B_k}{2} \cdot 2^2 \right) \\
    &=& 2B_k (X_k - I_k - R_k - 2B_k).
  \end{eqnarray*}
  The explanation is as follows.  $X_k$ is the sum of all $|C_j|^2$.
  Then, we break the sum over $j \not \in Z$ of $|C_j|^2$ into the
  cases when $C_j$ has type I, R, or B, in which $|C_j|$ is always 1,
  1, and 2, respectively.

  The total contribution from pairs of type BJ and JB is twice that
  from BJ alone, so it is $4B_k (X_k - I_k - R_k - 2B_k)$.

\item[Case RR.] Now $p_{ij} = 1$.  Any $|C_i|^2 |C_j|^2$ of this type
  is precisely $1^2 \cdot 1^2 = 1$.  The number of $C_i$ of type R is
  $R_k$, because the red matching consists of $\frac{R_k}{2}$ isolated
  red edges, which give $R_k$ isolated vertices in the blue color
  class.  So, the number of pairs $C_i, C_j$ of type RR with $i \neq
  j$ is $R_k \cdot (R_k - 1) = R_k^2 - O(n)$.  Thus the sum of
  $|C_i|^2 |C_j|^2$ is $R_k^2 - O(n)$.

\item[Cases RI, IR.] Again $p_{ij} = 1$.  Any $|C_i|^2 |C_j|^2$ of
  this type is precisely $1^2 \cdot 1^2 = 1$.  There are $R_k \cdot
  I_k$ pairs $C_i, C_j$ of type RI, and the same number of type IR, so
  the sum is $2 R_k I_k$.

\item[Cases RJ, JR.] Again $p_{ij} = 1$.  Let $Z$ be the set of indices
  $j$ such that $C_j$ has type J.  Since there are $R_k$ components
  $C_i$ of type R, the sum of $|C_i|^2 |C_j|^2$ over all pairs of type
  RJ is
  \begin{displaymath}
    R_k \sum_{j \in Z} 1^2 \cdot |C_j|^2
    \ = \ R_k (X_k - I_k - R_k - 2B_k),
  \end{displaymath}
  where we used the exact same calculation as in the case BJ for
  $\sum_{j \in Z} |C_j|^2$.  We double this to include the
  contribution from JR, and obtain a total sum of $2R_k (X_k - I_k -
  R_k - 2B_k)$.

\item[All other cases.] For all other pairs of types, our algorithm
  chooses a random color, so $p_{ij} = \frac{1}{2}$, and there is no
  difference between $S$ and $S'$.
\end{description}

Combining all of the above calculations, we express $\E{X_{k+1} - X_k
  \mid e_1, \ldots, e_k; \frac{1}{n} X_k < 210} = S$ in terms of $S'
\leq \big( \frac{X_k}{n} \big)^2 + O(n^{-1})$.
\begin{eqnarray*}
  S &=& S'
  \begin{array}[t]{l}
    - \frac{1}{n^2} \left[
      (4 B_k^2 - O(n))
      + 4 B_k I_k
      + 4B_k (X_k - I_k - R_k - 2B_k)
    \right] \\
    + \frac{1}{n^2} \left[
      (R_k^2 - O(n))
      + 2 R_k I_k
      + 2R_k (X_k - I_k - R_k - 2B_k)
    \right].
  \end{array} \\
  &=& S' +
  \frac{1}{n^2} \left[
    4B_k^2 - 4B_k X_k - R_k^2 + 2R_k X_k
  \right]
  + O(n^{-1}) \\
  &\leq& \left( \frac{X_k}{n} \right)^2 +
  \frac{1}{n^2} \left[
    4B_k^2 - 4B_k X_k - R_k^2 + 2R_k X_k
  \right]
  + O(n^{-1}),
\end{eqnarray*}
as desired.  \hfill $\Box$

\vspace{3mm}

By Lemma \ref{lem:track-partition}, $\frac{1}{n} B_k$ and $\frac{1}{n}
R_k$ track $b(t) = t e^{-4t}$, so Lemma \ref{lem:deltaX} indeed
indicates that the differential equation \eqref{xdiff} estimates
$\frac{1}{n} X_{tn}$.  We now prove this formally.  Our method uses
Hoeffding-Azuma, so we need bounded differences.  In our proof of
Theorem \ref{thm:evolve-suscep}, we achieved this by controlling the
distribution of the component sizes with the result of Spencer and
Wormald (Fact \ref{fact:SW}).

Recall that a graph has a $K, c$ component tail if for all positive
integers $s$, at most $Ke^{-cs}$-fraction of vertices lie in
components of order at least $s$.  In particular, the empty graph has
$K,c$ component tail with $K=e$ and $c=1$.  Fact \ref{fact:SW} then
ensures that after a period of random edge addition, the resulting
graph still has a $K', c'$ component tail.  However, the period only
lasts for about $0.5n$ edges when starting with the empty graph, and
our process needs to run for $1.06n$ rounds.  To work around this
issue, we use several iterations.

Define the sequence $t_0, \ldots, t_{19}$, by letting $t_0 = 0$, and
$t_{i+1} = t_i + \frac{1}{4x(t_i)}$, where $x(t)$ is the solution of
the differential equation \eqref{xdiff}.  The motivation for this
sequence is as follows.  Suppose we have already established that the
blue graph after $t_i n$ rounds has a $K_i, c_i$ component tail, and
its susceptibility $L$ is approximately $x(t_i)$, specifically, that
$L < 1.5x(t_i)$.  Then, we could apply Fact \ref{fact:SW} with $L =
1.5x(t_i)$, $K = K_i$, $c = c_i$, and $\gamma = \frac{1}{4}$, to to
conclude that after $t_{i+1} n$ rounds, even if all new edges were
colored blue, the blue graph would still have a $K_{i+1}, c_{i+1}$
component tail \whp.  This allows us to define sequences $K_0 \leq
\cdots \leq K_{19} = K'$ and $c_1 \geq \cdots \geq c_{19} = c'$.  We
confirmed numerically that $t_{19} > 1.06$, so this would allow us to
maintain a $K', c'$ component tail for $1.06n$ rounds.  Now we
formalize this heuristic, and prove our two-color avoidance theorem.

\vspace{3mm}


\noindent \textbf{Proof of Theorem \ref{thm:online-avoid-2}.}\, For
each $0 \leq k \leq 1.06 n$, let $\mathcal{E}_k$ be the event that all
of the following hold:
\begin{displaymath}
  \mathcal{E}_k
  \ = \
  \left\{
    \begin{array}{l}
      \big|\frac{1}{n} B_k -b\big(\frac{k}{n}\big)\big|
      \ \leq \ e^8 n^{-\frac{1}{3}}, \\
      \big|\frac{1}{n} R_k -b\big(\frac{k}{n}\big)\big|
      \ \leq \ e^8 n^{-\frac{1}{3}}, \\
      \frac{1}{n} X_k \ \leq \ x\big(\frac{k}{n}\big) + e^{\frac{500k}{n}} n^{-\frac{1}{4}}, \\
      \text{and the blue graph has a $K', c'$ component tail}.
    \end{array}
  \right.
\end{displaymath}
We define a supermartingale. Let
\begin{displaymath}
  Z_k \ = \ \begin{cases}
    X_k-nx\big(\frac{k}{n}\big)-e^{\frac{500k}{n}} n^{\frac{3}{4}}&\text{if }\cE_{k-1}\text{ occurs},\\
    Z_{k-1}&\text{otherwise}.
  \end{cases}
\end{displaymath}
We only consider $k\leq 1.06n$, and $x(t) \leq 209$ for all $t \leq
1.06$, so if $\mathcal{E}_k$ holds, we have $\frac{1}{n} X_k < 210$.
Then Lemma \ref{lem:deltaX} gives
\begin{align*}
&\E{Z_{k+1}-Z_k \mid e_1, \ldots, e_k, \mathcal{E}_k}\\
&\alignsp \leq \ \left( \frac{X_k}{n} \right)^2 +   \frac{1}{n^2} \left[4B_k^2 - 4B_k X_k - R_k^2 + 2R_k X_k\right]  + O(n^{-1})\\
&\alignsp \hspace{62pt}-n\left[x\left(\frac{k+1}{n}\right)-x\left(\frac{k}{n}\right)\right]
-\left[e^{\frac{500(k+1)}{n}} - e^{\frac{500k}{n}} \right]n^{\frac{3}{4}}.
\end{align*}
Now we estimate each term.  Since $\frac{X_k}{n} \leq
x\big(\frac{k}{n}\big) + e^{\frac{500k}{n}} n^{-\frac{1}{4}}$ and $k
\leq 1.06n$, we have $\big(\frac{X_k}{n}\big)^2 \leq
x^2\big(\frac{k}{n}\big) + 2 x\big(\frac{k}{n}\big)
e^{\frac{500k}{n}} n^{-\frac{1}{4}} + O(n^{-\frac{1}{2}})$.
Similarly, $\frac{B_k^2}{n^2}  = b^2\big(\frac{k}{n}\big) +
O(n^{-\frac{1}{3}})$, and the same estimate holds for
$\frac{R_k^2}{n^2}$.  Also, $\frac{1}{n} (2 R_k - 4 B_k) =
-2b\big(\frac{k}{n}\big) + O(n^{-\frac{1}{3}})$, so
\begin{displaymath}
  \frac{1}{n} (2 R_k - 4 B_k) \cdot \frac{X_k}{n}
  \ \leq \
  -2b\left(\frac{k}{n}\right) \left[
    x\left(\frac{k}{n}\right)
    - e^{\frac{500k}{n}} n^{-\frac{1}{4}}
  \right]
  + O(n^{-\frac{1}{3}}).
\end{displaymath}
From Taylor bounds similar to those in the proof of Lemma
\ref{lem:track-partition}, we have $x\big(\frac{k+1}{n}\big) -
x\big(\frac{k}{n}\big) = \frac{1}{n} x'\big(\frac{k}{n}\big) +
O(n^{-2})$ and $e^{\frac{500(k+1)}{n}} - e^{\frac{500k}{n}} =
\frac{500}{n} e^{\frac{500k}{n}} + O(n^{-2})$.  Combining all of
these bounds, and using $x' = x^2 + 3b^2 - 2bx$, the entire estimate
simplifies to
\begin{displaymath}
\E{Z_{k+1}-Z_k \mid e_1, \ldots, e_k, \mathcal{E}_k}
\ \leq \ \left[
  2x\left(\frac{k}{n}\right)
  + 2b\left(\frac{k}{n}\right)
  - 500
  \right] e^{\frac{500k}{n}} n^{-\frac{1}{4}}
  + O(n^{-\frac{1}{3}}),
\end{displaymath}
which is indeed less than zero for large $n$ because $b(t) = t
e^{-4t}$ is always less than 1, and $x(t) \leq 209$ for all $t \leq
1.06$.  Therefore $Z_0, \ldots, Z_{1.06 n}$ is a supermartingale.  Note
that $Z_0 = -n^{\frac{3}{4}}$.  Now because we are dealing with a
graph with a $K',c'$ tail, just as in the proof of Theorem
\ref{thm:evolve-suscep} we have $|Z_{k+1}-Z_k|=O(\log^2n)$ and then
the Hoeffding-Azuma inequality implies that for each $k \leq 1.06n$,
\begin{displaymath}
  \pr{ Z_k \geq 0 } \leq e^{-\Omega(n^{1/2} / \log^4 n)}.
\end{displaymath}
Therefore, by a union bound, \whp\ $Z_k < 0$ for all $k \leq 1.06 n$.
Also, Lemma \ref{lem:track-partition} implies that \whp,
$\big|\frac{B_k}{n}-b\big(\frac{k}{n}\big)\big| \leq e^8
n^{-\frac{1}{3}}$ and $\big|\frac{R_k}{n}-b\big(\frac{k}{n}\big)\big|
\leq e^8 n^{-\frac{1}{3}}$ for every $k \leq 1.06 n$.  Let
$\mathcal{E}$ be the conjunction of all of these high-probability
events.

To complete our argument, we show by induction that \whp, for each $0
\leq i \leq 19$, the blue graph after $t_i n$ rounds has a $K_i, c_i$
component tail.  The base case $i=0$ is trivial.  For the induction
step, suppose that it is true for $i$.  Condition on the blue graph
after $t_i n$ rounds having a $K_i, c_i$ component tail, as well as on
the event $\mathcal{E}$ that all $Z_k < 0$ and all $B_k$, $R_k$ are
concentrated.  Then, the same argument as in the conclusion of the
proof of Theorem \ref{thm:evolve-suscep} forces all $\mathcal{E}_k$ to
occur for $k \leq t_i n$, since $K_i \leq K'$ and $c_i \geq c'$.  In
particular, $\mathcal{E}_{t_i n}$ already implies that after $t_i n$
rounds, the blue graph has susceptibility $\frac{1}{n} X_{t_i n} \leq
x(t_i) + o(1) < 1.5 x(t_i)$.  Applying Fact \ref{fact:SW} with $L =
1.5 x(t_i)$, $K = K_i$, $c = c_i$, and $\gamma = \frac{1}{4}$, we see
that \whp, even if all new edges were colored blue, the blue graph
after $t_i n + \big(1 - \frac{1}{4}\big) \frac{n}{2 \cdot 1.5 x(t_i)}
= t_{i+1} n$ rounds would have a $K_{i+1}, c_{i+1}$ component tail.
This finishes the induction, so \whp\ the blue graph after $t_{19} n >
1.06 n$ rounds has a $K', c'$ component tail.  In particular, all
connected components are of order $O(\log n)$, so there is no giant in
the blue color class.  The same result follows for the red color class
by symmetry.  \hfill $\Box$

\section{Online creation of giants}
\label{sec:online-create}

Recall that the trivial bounds for the online creation of giants are
as follows.  No algorithm can create giants in all colors in fewer
than $(1-\epsilon)\frac{n}{2}$ total edges, because that is not even
enough to make a giant in the uncolored graph.  On the other hand,
if one randomly colors each incoming edge, then monochromatic giants
will appear after $(r + \epsilon) \frac{n}{2}$ total edges.  In this
section, we prove Theorems \ref{thm:online-embrace-lower} and
\ref{thm:online-embrace-upper}, which improve the above trivial
lower and upper bounds for the online creation of giants.

\subsection{Lower bound}

The previous argument iterated Fact \ref{fact:SW} to maintain the
component tail property, using a customized argument to control the
susceptibility for a specific algorithm.  In this section, we need
to consider an arbitrary coloring strategy, so we use our
general-purpose tool (Theorem \ref{thm:evolve-suscep}) to control
the susceptibility. This will establish a lower bound of $\Omega(n
\log r)$ for the number of edges required to create giants online in
each of $r$ color classes.  We need the following simple bound for
random graphs.

\begin{lemma}
  \label{lem:few-short-cycles}
  Let $\lambda$ be a constant.  The random graph $G_{n,p}$ with $p =
  \frac{\lambda}{n}$ contains at most $o\big( \frac{n}{\log n} \big)$
  cycles of length at most $\sqrt{\log n}$, \whp.
\end{lemma}

\noindent \textbf{Proof.}\, The expected number of cycles of length
$k$ in $G_{n,p}$ is at most $\frac{n^k}{2k} p^k =
\frac{\lambda^k}{2k}$, so the expected number of cycles of length at
most $\sqrt{\log n}$ is below $\sum_{k=3}^{\sqrt{\log n}}
\frac{\lambda^k}{2k}$.  If $\lambda \leq 1$, this is below $\sqrt{\log
  n}$.  Otherwise, it is below $\sqrt{\log n} \cdot
\lambda^{\sqrt{\log n}}$.  In both cases, the conclusion follows from
Markov's inequality.  \hfill $\Box$

\vspace{3mm}

Next, we need a worst-case bound on how large the susceptibilities of
different color classes can be when a graph is colored.

\begin{lemma}
  \label{lem:split-colors}
  Let $K, c$ be positive real constants.  Let $G$ be an $n$-vertex
  graph with a $K,c$ component tail.  Also assume that $G$ contains
  $o\big(\frac{n}{\log n}\big)$ cycles of length at most $\sqrt{\log
    n}$.  Consider any 2-coloring of the edges of $G$, and let
  $G^{(1)}$ and $G^{(2)}$ be the $n$-vertex subgraphs of $G$ obtained
  by keeping only edges in the first or second color, respectively.
  Then $S(G^{(1)}) + S(G^{(2)}) \leq S(G) + 1 + o(1)$.
\end{lemma}

\noindent \textbf{Proof.}\, Each component of $G^{(i)}$ is entirely
contained within a component of $G$, so we may break down the left
hand side by components of $G$.  Consider first the components of
$G$ which are larger than $\sqrt{\log n}$.  Since $G$ has a $K,c$
component tail, the number of vertices in such components is at most
$K e^{-c \sqrt{\log n}} n$.  The component tail also implies that
there is some constant $C$ such that all components of $G$ are
bounded by $C \log n$.  Since $S(G^{(1)}) + S(G^{(2)}) =
\frac{1}{n}\sum_v (C_v^{(1)} + C_v^{(2)})$, where $C_v^{(i)}$ is
number of vertices in the component of $G^{(i)}$ containing $v$, the
total contribution from vertices in components of $G$ with order at
least $\sqrt{\log n}$ is only $\frac{1}{n} \cdot K e^{-c \sqrt{\log
n}} n \cdot 2C\log n = o(1)$.

Next, consider the components of order at most $\sqrt{\log n}$ which
contain cycles.  Since the susceptibility is $\frac{1}{n}$ times the
sum of squares of component sizes, each component of this type
contributes at most $\frac{1}{n} \cdot 2(\sqrt{\log n})^2$ to
$S(G^{(1)}) + S(G^{(2)})$.  By assumption, $G$ only has
$o\big(\frac{n}{\log n}\big)$ cycles small enough to fit into these
components, so the number of such components is at most
$o\big(\frac{n}{\log n}\big)$.  Therefore, their total contribution to
$S(G^{(1)}) + S(G^{(2)})$ is at most $\frac{1}{n} \cdot
2(\sqrt{\log n})^2 \cdot o\big(\frac{n}{\log n}\big) = o(1)$.

The main contribution comes from the remaining components, which are
all trees.  Any tree $T$ in $G$ contributes $\frac{1}{n} \sum_{v \in
  T} |T|$ to $S(G)$.  We claim that it contributes at most
$\frac{1}{n} \sum_{v \in T} (|T| + 1)$ to $S(G^{(1)}) + S(G^{(2)})$,
i.e., the additional amount is at most $\frac{1}{n} |T|$.  Indeed,
$T$'s contribution to $S(G^{(i)})$ is precisely $\frac{1}{n}$ times
the sum of the sizes of the color-$i$ components that contain each
vertex $v \in T$.  Trees have the property that each pair of
vertices is connected by a unique path, so we can express the size
of the color-$i$ component containing $v$ as $\sum_{w \in T}
I_{v,w}^{(i)}$, where the indicator $I_{v,w}^{(i)}$ is 1 if the
unique path between $v$ and $w$ is monochromatic in color $i$, and 0
otherwise.  Hence, the total contribution of $T$ to $S(G^{(1)}) +
S(G^{(2)})$ is $\frac{1}{n} \sum_{v,w \in T} (I_{v,w}^{(1)} +
I_{v,w}^{(2)})$.  Since $T$ is a tree, the only time both indicators
$I_{v,w}^{(i)}$ can be 1 is when $w=v$.  So for each $v$, the sum
$\sum_{w \in T} (I_{v,w}^{(1)} + I_{v,w}^{(2)})$ is at most $|T| +
1$, as claimed.  Summing over all tree components, we see that their
total contribution to $S(G^{(1)}) + S(G^{(2)})$ exceeds $S(G)$ by at
most $\frac{1}{n}$ times the sum of the sizes of tree components,
which is at most 1.  Combining this with the contributions from
non-tree components above, we obtain $S(G^{(1)}) + S(G^{(2)}) \leq
S(G) + 1 + o(1)$, as desired.  \hfill $\Box$

\vspace{3mm}

Now we proceed to prove Theorem \ref{thm:online-embrace-lower},
using the previous two lemmas, and Theorem \ref{thm:evolve-suscep}
to control the evolution of susceptibility.  We will show that for
any $r$ which is a power of two, \whp\ no online algorithm can
create giants in all $r$ colors within $(c \log_2 r) n$ edges, where
$c \approx 0.043$.  This clearly implies the desired asymptotic
bound. Our calculated bound for $r=2$ will follow as a special case.

\vspace{3mm}

\noindent \textbf{Proof of Theorem \ref{thm:online-embrace-lower}.}\,
Let $C_0$ be the set of all $r = 2^t$ colors.  Let $\gamma$ be a
constant parameter which we will specify later.  The graph is
initially empty, with susceptibility $L_0 = 1$.  By Theorem
\ref{thm:evolve-suscep}, after $(1-\gamma) \frac{n}{2} L_0^{-1}$
edges, the graph formed by the union of those edges has $K_1, c_1$
component tail and susceptibility at most $\frac{L_0}{\gamma} + o(1)$
\whp.  Arbitrarily divide the colors into two groups of size $2^{t-1}$
each.  Lemmas \ref{lem:few-short-cycles} and \ref{lem:split-colors}
ensure that no matter how the edges were colored, one of the two color
groups determines a graph $G_1$ with susceptibility at most $L_1 +
o(1)$, where $L_1 = \frac{1}{2}\big(1 + \frac{L_0}{\gamma}\big)$.
Note that $G_1$ still has $K_1, c_1$-component tail, and let $C_1$ be
the set of $2^{t-1}$ colors we picked.

We iterate this procedure a total of $t$ times.  For example, in the
next step, we advance by $(1-\gamma) \frac{n}{2} L_1^{-1}$ more edges.
Even if all of them received colors in $C_1$ (i.e., were added to
$G_1$), the susceptibility of the graph determined by $C_1$-colors is
at most $\frac{L_1}{\gamma} + o(1)$ \whp, by Theorem
\ref{thm:evolve-suscep}.  Arbitrarily divide the colors of $C_1$ into
two groups of size $2^{t-2}$ each.  Again by Lemmas
\ref{lem:few-short-cycles} and \ref{lem:split-colors}, one of the two
color groups, say $C_2$, determines a graph $G_2$ with susceptibility
at most $L_2 + o(1)$, where $L_2 = \frac{1}{2}\big( 1 +
\frac{L_1}{\gamma} \big)$.

After $t$ iterations, we conclude that there is some single color $c$
such that the graph $G_t$ determined by all edges of color $c$ has a
$K_t, c_t$ component tail and susceptibility at most $L_t$.  A final
application of Theorem \ref{thm:evolve-suscep} implies that we can add
$\frac{n}{2} (L_t^{-1} - \epsilon)$ more random edges and still have
all components in color $c$ of order $O(\log n)$ \whp.

It remains to count the total number of edges which we have
accumulated.  The relationship between the $L_i$'s is $L_{i+1} =
\frac{1}{2}\big( 1 + \frac{L_i}{\gamma} \big) = \frac{1}{2} +
\frac{L_i}{2\gamma}$, so
\begin{eqnarray*}
  L_0 &=& 1, \\
  L_1 &=& \frac{1}{2} + \frac{1}{2\gamma}, \\
  L_2 &=& \frac{1}{2} + \frac{1}{4\gamma} + \frac{1}{4\gamma^2}, \\
  L_3 &=& \frac{1}{2} + \frac{1}{4\gamma} + \frac{1}{8\gamma^2} + \frac{1}{8\gamma^3}, \\
\end{eqnarray*}
and in general,
\begin{eqnarray*}
  L_t &=& \frac{1}{2} + \frac{1}{2(2\gamma)} + \frac{1}{2(2\gamma)^2} + \cdots
  + \frac{1}{2(2\gamma)^{t-1}} + \frac{1}{(2\gamma)^{t}} \\
  &<& 1 + \frac{1}{2\gamma} + \cdots + \frac{1}{(2\gamma)^{t}} \\
  &<& \left(1 - \frac{1}{2\gamma}\right)^{-1}.
\end{eqnarray*}
Thus, the total number of edges added (not even counting the final
step) is at least
\begin{displaymath}
  (1-\gamma)\frac{n}{2} \sum_{i=0}^{t-1} L_i^{-1}
  \ > \
  (1-\gamma)\frac{n}{2} \cdot t \left(1 - \frac{1}{2\gamma}\right).
\end{displaymath}
By routine calculus, the optimal choice for $\gamma$ is
$\frac{1}{\sqrt{2}}$, giving $(1-\gamma)\big(1 -
\frac{1}{2\gamma}\big) = \frac{3}{2} - \sqrt{2} \approx 0.086$.
Since $t = \log_2 r$, we indeed see that \whp, no online algorithm
can create giants in all colors within $0.043 n \log_2 r$ edges.
This completes the proof of the asymptotic bound.

For the specific case of $r=2$ colors, we can add the final batch of
$\frac{n}{2} (L_t^{-1} - \epsilon)$ random edges (here $t = 1$) to get
a specific bound which beats the trivial bound of $n/2$ edges.  Since
$L_1 = \frac{1}{2} \big( 1 + \frac{1}{\gamma} \big)$, this gives a
total edge count of
\begin{eqnarray*}
  (1-\gamma) \frac{n}{2} + \frac{n}{2} (L_1^{-1} - \epsilon)
  &=& \frac{n}{2} \left[ (1-\gamma) + \left(\frac{1}{2}\left(1 + \frac{1}{\gamma}\right)\right)^{-1} - \epsilon \right] \\
  &=& \frac{n}{2} \left[ (1-\gamma) + \frac{2\gamma}{\gamma + 1} - \epsilon \right].
\end{eqnarray*}
By routine calculus, the optimal choice for $\gamma$ is
$\sqrt{2}-1$. Therefore, \whp, no online algorithm can create giants
in both colors within $(2-\sqrt{2} - \epsilon)n$ edges, as claimed.
\hfill $\Box$

\subsection{Upper bound for many colors}

In this section, we present an online coloring algorithm which
creates giants in all $r$ color classes within roughly $\frac{n}{2}
\sqrt{r}$ edges.  The strategy is based on the classical fact that
there are infinitely many values of $r$ such that the edges of $K_r$
can be perfectly partitioned into cliques of order roughly
$\sqrt{r}$.

\begin{fact}
  \label{fact:partition-into-cliques}
  Let $r = q^2 + q + 1$ for some prime power $q$.  The edges of $K_r$
  can be partitioned into disjoint sets $E_1, \ldots, E_r$ such that
  each $E_i$ is precisely the edge set of some clique of order $q+1$.
\end{fact}

\noindent \textbf{Proof.}\, The projective plane of order $r = q^2 + q
+ 1$ is the finite geometry where points and lines correspond to
dimension-1 and dimension-2 subspaces of $\mathbb{F}_q^3$,
respectively.  This object contains exactly $\frac{q^3 - 1}{q-1} = q^2
+ q + 1$ points and the same number of lines, and has the property
that every pair of distinct points determines a unique line.

Identify the vertices of $K_r$ with the points of the projective
plane.  Let the $q+1$ vertices of the clique corresponding to $E_i$ be
the points contained in the $i$-th line of the projective plane.  The
edge partition property is then equivalent to the incidence property
of the projective plane.  \hfill $\Box$

\vspace{3mm}

We also need the giant component threshold in certain inhomogeneous
random graph models, where the edge probability is not uniformly $p$
at all $n \choose 2$ possible sites.  Instead, the probability of each
edge depends on the locations of its endpoints.  Bollob\'as, Janson,
and Riordan recently completed a far-reaching study of phase
transitions in these types of inhomogeneous models in \cite{BJR}.  We
use a special case of their work, regarding the specific model below.

Fix a symmetric $k \times k$ matrix $A = (a_{ij})$.  Let $G_{n,A}$ be
the $n$-vertex random graph defined as follows.  Split the $n$
vertices into $k$ groups of size $n/k$.  Between each pair of distinct
vertices, say from the $i$-th and $j$-th groups (where $i$ may equal
$j$), place an independent random edge with probability
$\frac{a_{ij}}{n}$.  Note that when $A = cJ_k$, where $J_k$ is the $k
\times k$ all-ones matrix, $G_{n,A}$ is the Erd\H{o}s-R\'enyi random
graph $G_{n,p}$ with $p = \frac{c}{n}$.

The following result was proved as Theorem 3.1 of \cite{BJR}.  Here,
the $L_2$ operator norm $\|B\|_2$ of a $k \times k$ matrix $B$ is
$\sup\{ \|Bx\|_2 : \|x\|_2 = 1 \}$, and the 2-norm of a vector $(x_1,
\ldots, x_k)$ is $\sqrt{\sum x_i^2}$.

\begin{fact}
  \label{fact:BJR}
  Let $A = (a_{ij})$ be a symmetric $k \times k$ matrix, and let
  $\overline{A}$ be its normalization $\big(\frac{a_{ij}}{k}\big)$.
  If $\|\overline{A}\|_2 > 1$, then $G_{n,A}$ contains a giant
  component \whp.
\end{fact}

\noindent \textbf{Remark 1.}\, In the same theorem, Bollob\'as,
Janson, and Riordan also proved the complementary result that when
$\|\overline{A}\|_2 \leq 1$, the largest component of $G_{n,A}$ is
$o(n)$ \whp.  However, we do not need this part for our analysis.

\vspace{3mm}

\noindent \textbf{Remark 2.}\, The $L_2$ operator norm of a real
symmetric matrix $A$ always equals its \emph{spectral radius}
$\rho(A)$, which is the maximum $|\lambda_i|$ over all eigenvalues
$\lambda_i$.  Indeed, $A$ is diagonalizable with an orthonormal basis
of real eigenvectors, so let the eigenvalues and eigenvectors be
$\lambda_1, \ldots, \lambda_k$ and $v_1, \ldots, v_k$, respectively.
Expressing any vector $x$ in this basis as $\sum c_i v_i$, we have
that the condition $\|x\|_2 = 1$ is precisely $\sum c_i^2 = 1$, and
$\|A x\|_2 = \sqrt{\sum \lambda_i^2 c_i^2}$.  Therefore, $\|Ax\|_2$
has maximum value equal to the largest absolute value of an
eigenvalue.

\vspace{3mm}

\noindent \textbf{Remark 3.}\, The Erd\H{o}s-R\'enyi model $G_{n,p}$
with $p = \frac{c}{n}$ corresponds to $G_{n,A}$ with $A = c J_k$.
The normalized matrix $\overline{A} = \frac{c}{k} J_k$ has
eigenvalues $c$ and 0, so Fact \ref{fact:BJR} implies the classical
result that the giant component appears after $p = \frac{1}{n}$.

\vspace{3mm}

We use this to study the $k$-partite random graph $G_{n,p}^{(k)}$,
which has $n$ vertices split into equal groups of size $\frac{n}{k}$,
and independent random edges with probability $p = \frac{c}{n}$
between pairs of vertices from distinct groups.  In the above
framework, this is $G_{n,A}$ with $A = c(J_k - I_k)$.

\begin{corollary}
  \label{cor:giant-punctured}
  Let $k \geq 2$ be a positive integer, and let $c > \frac{k}{k-1}$ be
  a real number.  Then the $k$-partite random graph $G_{n,p}^{(k)}$
  with $p = \frac{c}{n}$ contains a giant component \whp.
\end{corollary}

\noindent \textbf{Proof.}\, By Fact \ref{fact:BJR} and our second
remark, the problem reduces to determining the eigenvalues of
$\overline{A} = \frac{c}{k}(J_k - I_k)$.  These are precisely
$\frac{c}{k}(k-1)$ and $\frac{c}{k}(0-1)$, so since $k \geq 2$, the
giant component appears once $c > \frac{k}{k-1}$.  \hfill $\Box$

\vspace{3mm}

We are now ready to state our algorithm and prove its effectiveness.
Note that a coloring algorithm that produces giants in $r$ colors
trivially gives coloring algorithms for any $r' < r$ as well, simply
by using the first color whenever any color beyond $r'$ was to be
used.  So, Theorem \ref{thm:online-embrace-upper} is a consequence of
the following more precise formulation, combined with the Prime Number
Theorem and Fact \ref{fact:coupling}.

\begin{customtheorem}
  Let $r = q^2 + q + 1$ for some prime power $q$.  There is an online
  algorithm such that for any $\epsilon > 0$, \whp\ all $r$ color classes
  contain giant components within $\big(\frac{r}{q} + \epsilon\big)
  \frac{n}{2}$ edges.
\end{customtheorem}

\noindent \textbf{Proof.}\, Arbitrarily partition the $n$ vertices
into $r$ sets $V_1, \ldots, V_r$, each of size $\frac{n}{r}$.  By
Fact \ref{fact:partition-into-cliques}, there is a partition $E_1
\cup \ldots \cup E_r$ of the edges of $K_r$, such that each $E_t$ is
precisely the edge set of some clique of order $q+1$.  Our online
coloring algorithm is then as follows.  Usually, the incoming edge
will have endpoints in distinct parts $V_i$ and $V_j$.  In that
case, color the edge with the index $t$ of the $E_t$ which contains
the edge $ij$ in the partitioned graph $K_r$.  Otherwise, if the
incoming edge is spanned by a single $V_i$, then discard the edge
entirely.  Note that this is even stronger than coloring it, because
we will now find giants without using those edges at all.

Our algorithm disregards the entire history of the process, since the
color of each edge is a function of the locations of its endpoints.
In particular, the order of the edges is irrelevant, so the
performance only depends on the final edge set.  Thus, by Fact
\ref{fact:coupling}, it suffices to show that if this strategy is
applied to $G_{n,p}$ with $p = \big(\frac{r}{q} +
\frac{\epsilon}{2}\big) \frac{1}{n}$, then it creates giants in all
colors \whp.  By passing to this independent model, each color class
itself becomes a $(q+1)$-partite random graph $G_{n', p}^{(q+1)}$, on
only $n' = \frac{n}{r} (q+1) \approx \frac{n}{\sqrt{r}}$ vertices.
Indeed, $E_t$ is the edge set of a clique on some set $S$ of $q+1$
vertices of $K_r$, so the edges that receive color $t$ are precisely
those with endpoints in some $V_i$ and $V_j$ with $i \neq j$ and $i,j
\in S$.

Finally, we can apply Corollary \ref{cor:giant-punctured} with $k =
q+1$, since $p = \frac{c'}{n'}$ with $c' = \big(\frac{r}{q} +
\frac{\epsilon}{2}\big) \frac{1}{n} \cdot \frac{n}{r} (q+1) >
\frac{q+1}{q}$.  Therefore, each individual color class contains a
giant component \whp.  Taking a union bound over all $r$ (finitely
many) color classes finishes the proof.  \hfill $\Box$

\subsection{Upper bound for 2 colors}

To adapt our strategy from the previous section to the case $r=2$, we
must specify symmetric 0-1 matrices $A_1$ and $A_2$ which sum to the
$k \times k$ all-ones matrix $J_k$.  We then split the vertices into
$k$ equal parts $V_1, \ldots, V_k$, and color an edge with endpoints
in some $V_i, V_j$ with color 1 if the $ij$-entry of $A_1$ is 1, and
color 2 otherwise.

Then, after applying this strategy to the edges of $G_{n,p}$ with $p =
\frac{c}{n}$, the $i$-th color class is a copy of $G_{n, cA_i}$.  By
the second remark after Fact \ref{fact:BJR}, this contains a giant
component when the spectral radius $\rho(\frac{c}{k} A_i)$ exceeds 1.
Since our objective is to create giants in both colors as rapidly as
possible, we want to select $A_1$ and $A_2$ such that $A_1 + A_2 =
J_k$, but $\min\{\rho(A_1), \rho(A_2)\}$ is as large as possible.
This appears to be a nontrivial problem, but one simple way to choose
the matrices is to let $A_1$ have 1's in the top-left $t \times t$
submatrix, and 0's everywhere else.  This leads to the following
bound.

\begin{proposition}
  \label{prop:0.75create}
  For every $\epsilon>0$, it is possible to create giants in two colors online within $\big(
  \frac{3}{4} + \epsilon \big) n$ rounds \whp.
\end{proposition}

\noindent \textbf{Proof sketch.}\, Since $A_1$ is just $J_t$ embedded
in an all-zeros matrix, its spectral radius is precisely $t$.  Next,
note that $A_2 = J_k - A_1$ has rank 2, so it has at most 2 nonzero
eigenvalues $\lambda_1, \lambda_2$.  The trace of $A_2$ is $k-t$, so
$\lambda_1 + \lambda_2 = k-t$.  Also, the main diagonal of $A_2^2$ has
its first $t$ entries equal to $k-t$, and the remaining $k-t$ entries
equal to $k$, giving $\text{tr}(A_2^2) = t(k-t) + (k-t)k = k^2 - t^2$.
This trace also equals $\lambda_1^2 + \lambda_2^2$, because the
nonzero eigenvalues of $A_2^2$ are $\lambda_1^2$ and $\lambda_2^2$.
Solving this system of equations, one finds that the largest
eigenvalue of $A_2$ is $\frac{1}{2} ( k-t + \sqrt{k^2 + 2kt -3t^2} )$.
Recall that the largest eigenvalue of $A_1$ is $t$, and we wanted the
largest possible $\min\{\rho(A_1), \rho(A_2)\}$.  Routine calculus
shows that the optimal choice of $t$ is $\frac{2}{3}k$, giving both
$\rho(A_i) = \frac{2}{3}k$.  So, we choose the particular $3 \times 3$
matrices
\begin{displaymath}
  A_1 =
  \left(
    \begin{array}{ccc}
      1 & 1 & 0 \\
      1 & 1 & 0 \\
      0 & 0 & 0
    \end{array}
  \right),
  \quad \quad \quad
  A_2 =
  \left(
    \begin{array}{ccc}
      0 & 0 & 1 \\
      0 & 0 & 1 \\
      1 & 1 & 1
    \end{array}
  \right).
\end{displaymath}
Therefore, as we remarked at the beginning, Fact \ref{fact:BJR} shows
that when this strategy is applied to $G_{n,p}$ with $p =
\frac{c}{n}$, both colors will contain giant components if their
spectral radii $\rho(\frac{c}{k} A_i)$ exceed 1, i.e., once $c >
\frac{3}{2}$.  By Fact \ref{fact:coupling}, this happens after $\big(
\frac{3}{2} + \epsilon \big) \frac{n}{2}$ rounds, so we are done.
\hfill $\Box$

\vspace{3mm}

\noindent \textbf{Remark.}\, Although the partition we chose may
appear na\"ive, there is evidence to suggest that it may be optimal.
Note that if we ignore the main diagonal (an effect that can be made
negligible by choosing large $k$) and seek $A_1 + A_2 = J_k - I_k$,
then $A_1$ and $A_2$ are the adjacency matrices of a graph and its
complement.

Several researchers have studied the question of bounding the sum of
the spectral radii of the adjacency matrices of complementary graphs
(see \cite{HS, Li, Nik, Nik2, Nos, Zhou}).  In particular, Nikiforov
recently conjectured in \cite{Nik} that the sum of these two
spectral radii is always at most $\frac{4}{3} k + O(1)$, where $k$
is the number of vertices. If true, this would imply that
$\min\{\rho(A_1), \rho(A_2)\} \leq \frac{2}{3} k + O(1)$, which our
construction achieved.  In fact, in his extremal example, one graph
was a clique on a subset of the vertices, which is essentially the
same as our construction. So, perhaps $\frac{3}{4}n$ is the limit of
what can be achieved by any strategy as above.

\vspace{3mm}

Next, we prove Theorem \ref{thm:online-embrace-upper-2}, which shows
that by making the strategy more adaptive, one can create giants even
faster.  The algorithm in the proof of Proposition
\ref{prop:0.75create} fixed a subset $R$ of vertices in advance, and
used the first color whenever an edge was spanned by $R$.  The key
idea is to let the subset $R$ depend on the outcomes of the first few
rounds.  To analyze this strategy, we will need two results of Spencer
and Wormald.  These are Theorems 1.1 and 3.1 of \cite{SW}, restated in
equivalent form using Fact \ref{fact:coupling}.

\begin{fact}
  \label{fact:SW-suscep-lower}
  Let $0 < t < \frac{1}{2}$ be a fixed parameter.  Then there exist
  constants $K,c$ such that \whp, the graph on $n$ vertices formed by
  $t n$ independent random edges has susceptibility $\frac{1}{1-2t} +
  o(1)$, and a $K,c$ component tail.
\end{fact}

\begin{fact}
  \label{fact:SW-blow-up}
  Let $L, K, c, \epsilon$ be positive real numbers.  Let $G$ be a
  graph on $n$ vertices with a $K,c$ component tail and $S(G) = L$.
  Then, after adding $(1+\epsilon) \frac{n}{2L}$ more independent
  random edges, the resulting graph contains a giant component \whp.
\end{fact}

\vspace{3mm}

\noindent \textbf{Proof of Theorem
  \ref{thm:online-embrace-upper-2}.}\, Let the colors be red and blue.
We state the coloring strategy in terms of a constant parameter $t$,
which we can optimize at the end.  (The best choice turns out to be $t
\approx 0.189$.)  For the first $tn$ rounds, color all edges red.
Then, permanently fix $R$ to be the set of all vertices incident to a
red edge at that time.  Color each future edge red whenever both
endpoints lie in $R$, and blue otherwise.

Let $\alpha = \frac{|R|}{n}$.  The argument at the beginning of the
proof of Lemma \ref{lem:track-partition} shows that $\alpha = (1 -
e^{-2t} + o(1))$ \whp.  Let us analyze how many rounds are required
for a red giant to appear.  By Fact \ref{fact:SW-suscep-lower}, the
(completely red) graph $G$ at time $tn$ has susceptibility $S(G) =
\frac{1}{1-2t} + o(1)$ \whp, so the sum of the squares of its
component sizes is $\big(\frac{1}{1-2t} + o(1)\big)n$.  Let $G_R$ be
the subgraph of $G$ induced by $R$.  The sum of the squares of the
components in $G_R$ is precisely $S(G)n - (1 - \alpha)n$, because all
components of $G$ outside $R$ are singletons.  Therefore, since $G_R$
has $\alpha n$ vertices, its susceptibility $L$ is:
\begin{displaymath}
  L
  \ = \
  \frac{1}{\alpha n} [S(G)n - (1-\alpha)n]
  \ = \
  (1+o(1)) \frac{1}{\alpha} \left[ \frac{1}{1-2t} - e^{-2t} \right].
\end{displaymath}
Then, by Fact \ref{fact:SW-blow-up},
\whp\ the red graph will contain a giant component after $(1+\epsilon)
\frac{|R|}{2L}$ more random edges are added with both endpoints in
$R$.  By a standard coupling as in Fact \ref{fact:coupling}, this
happens after $(1+\epsilon) \frac{|R|}{2L} \cdot \alpha^{-2}$ more
rounds \whp, since each incoming edge falls within $R$ with
probability $\alpha^2$.  Substituting $|R| = \alpha n$, we find that a
red giant appears after a grand total of $tn + (\frac{1}{2 \alpha L} +
\epsilon) n = \big[ t + \frac{1}{2}\big( \frac{1}{1-2t} - e^{-2t}
\big)^{-1} + \epsilon \big] n$ rounds \whp.

To analyze the blue graph, observe that by a similar coupling to Fact
\ref{fact:coupling}, after $tn + (1+\epsilon) \frac{cn}{2}$ rounds the
blue graph contains $G_{n, cA}$ \whp, where $A$ is the $n \times n$
matrix with 0's in the top-left $|R| \times |R|$ submatrix, and 1's
everywhere else.  Plugging $|R| = \alpha n$ into the eigenvalue
calculation from the proof of Proposition \ref{prop:0.75create}, we
see that the largest eigenvalue of $A$ is $\frac{n}{2} \big( 1 -
\alpha + \sqrt{1 + 2\alpha - 3\alpha^2} \big)$.  Thus, Fact
\ref{fact:BJR} implies that \whp, the giant component appears in the
blue graph once $c$ surpasses $\frac{2}{1 - \alpha + \sqrt{1 + 2\alpha
    - 3\alpha^2}} + \epsilon$, i.e., when the total number of rounds
exceeds $tn + \frac{1+\epsilon}{1 - \alpha + \sqrt{1 + 2\alpha -
    3\alpha^2}}n$.

Since $\alpha = 1 - e^{-2t} + o(1)$, it is now routine to numerically
optimize $t$.  It turns out that the best choice is $t \approx 0.189$,
which gives $\alpha \approx 0.314$.  Then, both of the bounds at the
ends of the previous two paragraphs are satisfied after $0.733 n$
rounds, completing the proof.  \hfill $\Box$

\section{Concluding remarks}
In this paper we have introduced several rather natural algorithmic
variants of the classical problem of the appearance of the giant
component in a random graph/process. As expected, the offline cases of
these problems appear to be much more accessible, and indeed we
managed to solve both the avoidance and the embracing versions
asymptotically for any fixed $r$. The online case seems to be more
challenging; there we showed that in all cases one can do better than
the trivial algorithms that randomly color each incoming edge, but for
creating giants, rather sizable gaps remain.

It would certainly be nice to settle the case of two colors for
creating and avoiding giants online in both color classes, but that
could be difficult.  A more approachable problem might be to close the
asymptotic gap between the lower bound of $\Omega(\log r)\cdot n$ and
the upper bound of $O(\sqrt{r})\cdot n$ for the question of creating
giants in $r$ colors.  In particular, can one show a lower bound of
the form $r^an$ for some positive constant $a$?

Another, perhaps more technical, issue that we would like to see
settled is the nature of an algorithm for avoiding giants online.  Our
online avoidance algorithm is randomized. Is there a deterministic
strategy that matches its performance in the online setting?




\end{document}